\documentclass{article}
\usepackage{geometry}
\usepackage{amsmath, amssymb}
\usepackage{xcolor}
\usepackage{amsthm}
\usepackage[square,sort,comma,numbers]{natbib}
\usepackage{hyperref}
\hypersetup{
    colorlinks=true,
    linkcolor=blue,
    citecolor=blue
}
\usepackage{graphicx}
\usepackage[utf8]{inputenc}
\usepackage[english]{babel}

\newtheorem{theorem}{Theorem}

\newtheorem{definition}{Definition}
\newtheorem{proposition}{Proposition}
\newtheorem{corollary}{Corollary}
\newtheorem{lemma}{Lemma}
\newtheorem{remark}{Remark}
\newcommand{\indep}{\perp \!\!\! \perp}

\newcommand{\numberthis}{\addtocounter{equation}{1}\tag{\theequation}}

\title{Rate of convergence of the smoothed empirical Wasserstein distance}
\author{Adam Block, Zeyu Jia, Yury Polyanskiy, Alexander Rakhlin\\Massachusetts Institute of Technology}

\def\eqdef{\triangleq}
\def\EE{\mathbb{E}}
\def\PP{\mathbb{P}}
\def\dperp{\perp\!\!\!\perp}
\def\QQ{\mathbb{Q}}
\def\mreals{\mathbb{R}}
\long\def\nbyp#1{\textcolor{red}{[{\bf YP:} #1]}}

\def\my{\mathbf{y}}

\def\calN{\mathcal{N}}

\begin{document}
\maketitle

\begin{abstract}
	Consider an empirical measure $\PP_n$ induced by $n$ iid samples from a $d$-dimensional $K$-subgaussian distribution $\PP$ and let $\gamma = \mathcal{N}(0,\sigma^2 I_d)$ be the isotropic Gaussian measure. We study the speed of convergence of the smoothed Wasserstein distance $W_2(\PP_n * \gamma, \PP*\gamma) = n^{-\alpha + o(1)}$ with $*$ being the convolution of measures. For $K<\sigma$ and in any dimension $d\ge 1$ we show that $\alpha = {1\over2}$. For $K>\sigma$ in dimension $d=1$ we show that the rate is slower and is given by $\alpha = {(\sigma^2 + K^2)^2\over 4 (\sigma^4 + K^4)} < 1/2$. This resolves several open problems in~\cite{goldfeld2020convergence}, and in particular
	precisely identifies the amount of smoothing $\sigma$ needed to obtain a parametric rate. 
	In addition, for any $d$-dimensional $K$-subgaussian distribution $\mathbb{P}$, we also establish that $D_{KL}(\PP_n * \gamma \|\PP*\gamma)$ has rate
	$O(1/n)$ for $K<\sigma$ but only slows
	down to $O({(\log n)^{d+1}\over n})$ for $K>\sigma$. The surprising difference of the
	behavior of $W_2^2$ and KL implies the failure of $T_{2}$-transportation inequality when
	$\sigma < K$. Consequently, it follows that for $K>\sigma$ and $d=1$  the log-Sobolev inequality (LSI) for the Gaussian mixture $\mathbb{P} *
	\gamma$ cannot hold. This closes an open problem in~\cite{wang2016functional}, who
	established the LSI under the condition $K<\sigma$ and asked if their bound can be
	improved.
\end{abstract}

\section{Introduction and main results}


Given $n$ iid samples $X_1,\ldots, X_n$ from a probability measure $\PP$ on $\mreals^d$ let us
denote by $\PP_n = {1\over n} \sum_{i=1}^n \delta_{X_i}$ the  
empirical distribution. As $n\to \infty$ it is well known that $\PP_n \to \PP$ according to many different notions of convergence.
The literature on the topic is voluminous. Here we are interested in convergence in
Wasserstein $W_p$-distances, cf.~\cite[Chapter 1]{villani2003topics}, defined for $p\ge 1$ as
$$ W_p(\PP,\QQ)^p = \inf_{P_{X,Y}}\{ \EE[\|X-Y\|^p]: P_X = \PP, P_Y = \QQ\}\,,$$
where $\|\cdot\|$ is Euclidean norm. Already in~\cite{dudley1969speed} it was shown that 
$$ W_1(\PP_n, \PP) = \Theta(n^{-1/d})\,,$$
for $d\ge 3$ and compactly supported $\PP$ absolutely continuous with respect to Lebesgue measure.
Dudley's technique relied on the characterization (special to $p=1$) of $W_1$ as the supremum over expectations of Lipschitz functions. His idea of recursive partitioning was cleverly adapted to the
realm of couplings in~\cite{boissard2014mean}, recovering Dudley's convergence rate of $n^{-1/d}$ also for $p>1$. For standard Gaussian distributions, \cite{bobkov2019one} settles the rate of $\Theta(\frac{\log\log n}{n})$ for one-dimensional distributions, \cite{berthet2020exact} proves that the constant in the $\Theta$ is $1$, i.e. $\frac{\EE[W_1(\PP_n, \PP)]}{(\log\log n)/n} = 1 + o(1)$. \citep{ledoux2019optimal, talagrand2018scaling} settles the rate of $\Theta(\frac{\log^2 n}{n})$ for two-dimensional distributions, and \cite{bobkov2019one} settles the rate of $\Theta(\frac{1}{n(\log n)^{p/2}})$ for $d$ dimensional distributions with $d\ge 3$. See~\cite{dereich2013constructive,fournier2015rate,weed2019sharp} for more on this line of work, and also for a thorough survey of the recent literature. 

Somewhat surprisingly, it was discovered in~\cite{goldfeld2020convergence} that the rate of convergence improves all the way to (dimension-independent)
$n^{-1/2}$ if one merely regularizes both $\PP_n$ and $\PP$ by convolving with the Gaussian
density.\footnote{Of course, the price to pay for this fast rate is a constant in front of
$n^{-1/2}$, which can be exponential in $d$ for certain $\PP$, cf~\cite{goldfeld2020convergence}.} 
More precisely, let $\varphi_{\sigma^{2}}(x) \eqdef (2\pi \sigma^2)^{-d/2}
e^{-{\frac{\|x\|^2}{2\sigma^2}}}$ be the density of $\mathcal{N}(0,\sigma^2 I_d)$, and for any
probability measure $\PP$ on $\mreals^d$ we define the convolved measure via
  \begin{equation*}
       \PP * \mathcal{N}(0,\sigma^2 I_d) (E) = \int_{E} dz \EE\left[\varphi_{\sigma^{2}}(X-z)\right], \quad X\sim\mathbb{P},
  \end{equation*}
  where $E$ is any Borel set. Then~\cite[Prop. 6]{goldfeld2020convergence} shows
\begin{equation}\label{eq:w2_orig}
  	\EE[W_2^2(\PP_n * \mathcal{N}(0,\sigma^2 I_d), \PP * \mathcal{N}(0,\sigma^2 I_d))] \le {C(d,\sigma,K)\over n}\,,
\end{equation}  
  whenever $\PP$ is $K$-subgaussian and $K<{\sigma\over2}$. We recall that $X\sim\PP$ is
  $K$-subgaussian if
  \begin{equation}\label{eq:subgauss}
  	\EE[e^{(\lambda, X-\EE[X])}] \le e^{{1\over2} K^2\|\lambda\|^2} \qquad \forall \lambda \in
  \mreals^d\,.
\end{equation}  
  Note that in~\eqref{eq:w2_orig} constant $C$ does not depend on $\PP$. Estimate~\eqref{eq:w2_orig} is
  most exciting for large $d$, but even for $d=1$ and $\PP=\mathcal{N}(0,1)$ it is non-trivial as
  $\EE[W_2^2(\PP_n,\PP)] = \Theta({\log\log n\over n})$ (see \citep{bobkov2019one,berthet2020exact}). Another surprising feature
  is~\cite[Corollary 2]{goldfeld2020convergence}: for $K \ge \sqrt{2}\sigma$ there exists a
  $K$-subgaussian distribution $\PP$ in  $\mreals^1$ such that 
  	\begin{equation}\label{eq:w2_orig_lb}
		\lim_{n\to\infty} n\EE[W_2^2(\PP_n * \mathcal{N}(0,\sigma^2 I_d), \PP *\mathcal{N}(0,\sigma^2 I_d))]=\infty,
\end{equation}
  where the expectation is with respect to $n$ samples according to $\mathbb{P}$. We say that the rate of convergence is ``parametric'' if~\eqref{eq:w2_orig} holds and otherwise
  call it ``non-parametric''. Thus, the results of~\cite{goldfeld2020convergence} show that
  parametric rate for smoothed-$W_2$ is only attained by sufficiently light-tailed distributions
  $\PP$, e.g. the subgaussian constant of distribution $\PP$ is less than the scale of noise over two. 

In this paper we prove three principal results:
\begin{enumerate}
\item Theorem~\ref{th:w2_threshold} resolves the gap between the location of the parametric and
non-parametric region: it turns out that for $K<\sigma$ we always have~\eqref{eq:w2_orig}, while
for $K>\sigma$ we have~\eqref{eq:w2_orig_lb} for some $K$-subgaussian distribution $\PP$ in $\mreals^1$. 
  (We remark that for $W_1$ we always have parametric rate $n^{-1/2}$ for all
  $K,\sigma>0$, cf~\cite[Proposition 1]{goldfeld2020convergence}.)
\item In the region of non-parametric rates ($K>\sigma$) a natural question arises: what
rates of convergence are possible? In other words, what is the value of
\begin{equation}\label{eq:rho_def}
	\alpha = \alpha(K,\sigma, d) \eqdef  \inf_{\PP-K\text{-subgaussian}}\lim_{n\to\infty} - {\log \EE[W_2(\PP_n
*\mathcal{N}(0, \sigma^{2}I_{d}), \PP*\mathcal{N}(0, \sigma^{2}I_{d}))] \over \log n}. 
\end{equation}
Previously, it was only known that ${1\over 4} \le \alpha \le {1\over 2}$ for all $K>\sigma$ (note
that~\eqref{eq:w2_orig_lb} strongly suggests but does not formally imply $\alpha < {1\over 2}$). The lower bound $1/4$ of $\alpha$ follows from \cite{goldfeld2020convergence}. Theorem~\ref{thm-1D} shows that for $d=1$ we have 
$$  \alpha(K,\sigma,d=1) = {(\sigma^{2} + K^{2})^{2} \over 4(\sigma^{4} + K^{4}) }\,. $$

\item We can see that for a class of $K$-subgaussian distributions, the convergence rate of
$W_2(\PP_n*\mathcal{N}(0, \sigma^{2}I_{d}),\PP*\mathcal{N}(0, \sigma^{2}I_{d}))$ changes from $n^{-1/4}$ to $n^{-1/2}$ as
$\sigma$ increases from $0$ to $K$, after which the rate remains $n^{-1/2}$. Our final result
(Theorem~\ref{thm-KL})
shows that, despite being intimately related to $W_2$, the Kullback-Leibler (KL) divergence behaves
rather differently: For all $K$-subgaussian $\PP$ we have 
	\begin{equation}\label{eq:kl_rate}
		\EE[D_{KL}(\PP_n*\mathcal{N}(0, \sigma^{2}I_{d}) \| \PP *\mathcal{N}(0, \sigma^{2}I_{d}))] \le 
		\begin{cases} {C(\sigma,K,d) \log^{d}n\over n}, & K > \sigma\\
				{C(\sigma,K,d)\over n}, & K < \sigma
		\end{cases}
	\,, 
	\end{equation}	
	where $D_{KL}(\mu\|\nu) = \int d\nu f(x) \log f(x), f \eqdef {d\mu \over d\nu}$ whenever $\mu$
	is absolutely continuous with respect to $\nu$. Now from the proof of
	Theorem~\ref{th:w2_threshold} we also know that for $K>\sigma$, KL-divergence is $\omega
	({1\over n})$. Thus, while at $K>\sigma$ both $W_2$ and KL switch to the non-parametric regime,
	the $W_2$ distance experiences a polynomial slow-down in rate, while KL only gets hit by
	(at most) a poly-logarithmic penalty.
\end{enumerate}
	
To better understand the relationship between the $W_2$ results and the KL one, let us recall an
important result of Talagrand (known as $T_2$-transportation inequality). A probability measure
$\nu$ is said to satisfy the $T_2$ inequality if there exists a finite constant $C$ such that 
$$ \forall \QQ: \quad W_2^2(\QQ, \nu) \le C\cdot D_{KL}(\QQ\|\nu)\,.$$
The infimum over all such constants is denoted by $T_2(\nu)$. Talagrand originally demonstrated that $T_2(\mathcal{N}(0, \sigma^{2}I_{d}))<\infty$    
It turns out that $T_2(\PP*\mathcal{N}(0, \sigma^{2}I_{d}))<\infty$ as well for compactly supported $\PP$~\cite{zimmermann2013logarithmic} and $K$-subgaussian $\PP$ with $K<\sigma$~\cite{wang2016functional} (in fact, both papers establish
a stronger log-Sobolev inequality (LSI)). Also in~\citep{chen2021dimension}, sharper LSI constants for distribution $\PP * \mathcal{N}(0, \sigma^{2}I_{d})$ are given where $\PP$ is with compact support or subgaussianity.  

Now comparing~\eqref{eq:kl_rate} and the lower bound for all $K>\sigma$ established
in Theorem~\ref{thm-1D} we discover the following.

\begin{corollary}\label{cor:lsi} For any $K>\sigma$ there exists a $K$-subgaussian $\PP$ on $\mreals^1$ such that
$\PP*\mathcal{N}(0, \sigma^{2})$ does not satisfy $T_2$-transportation inequality (and hence does not satisfy
the LSI either), that is $T_2(\PP*\mathcal{N}(0, \sigma^{2}))=\infty$.
\end{corollary}
We remark that it is straightforward to show that
	$$ \sup \{T_2(\PP*\mathcal{N}(0, \sigma^{2})): \PP \mbox{~--~$K$-subgaussian}\} = \infty $$
by simply considering $\PP=(1-\epsilon)\delta_{0} + \epsilon \delta_{N}$ for $\epsilon\to0$ and
$N\to\infty$ (cf. Appendix~\ref{apx:t2}). However, each of these measures has $T_2<\infty$.
Evidently, our corollary proves a stronger claim. 

Incidentally, this strengthening resolves an open question stated in~\cite{wang2016functional},
who proved the LSI (and $T_2$)
for $\PP*\mathcal{N}(0, \sigma^{2})$ assuming $\EE[e^{aX^2}]<\infty$ holds for some $a>{1\over 2\sigma^2}$, where $X\sim\PP$. They raised a question whether this threshold can be reduced, and our Corollary shows
the answer is negative. Indeed, one only needs to notice that whenever $X\sim \PP$ is $K$-subgaussian it satisfies 
	\begin{equation}\label{eq:ww_co}
		\EE\left[e^{a X^2}\right] < \infty \quad \forall a < {1\over 2K^2}\,,
\end{equation}
which is proved in \cite[p. 26]{boucheron2013concentration}.

\subsection{Main results and proof ideas}

  Our first result is the following:
  \begin{theorem}\label{th:w2_threshold}
      If $K<\sigma$, then for any $K$-subgaussian distribution $\PP$, we have 
      $$\mathbb{E}\left[W_2^2(\PP_n*\mathcal{N}(0, \sigma^{2}I_{d}), \PP*\mathcal{N}(0, \sigma^{2}I_{d}))\right] = \mathcal{O}\left(\frac{1}{n}\right),$$
      where $\PP_n$ is the empirical measure of $\PP$ with $n$ samples, and the expectation is over these $n$ samples.
  If $K > \sigma$, then there exists a $K$-subgaussian distribution $\PP$ such that 
  $$\mathbb{E}\left[W_2^2(\PP_n*\mathcal{N}(0, \sigma^2I_d), \PP*\mathcal{N}(0, \sigma^2I_d))\right] = \omega\left(\frac{1}{n}\right).$$
  \end{theorem}
  
  \paragraph{Previous results.}~\cite{goldfeld2020convergence} shows when $K<\sigma/2$,
  $\mathbb{E}\left[W_2^2(\PP_n*\mathcal{N}(0, \sigma^2I_d), \PP*\mathcal{N}(0, \sigma^2I_d))\right]$
  converges with rate $\mathcal{O}\left(\frac{1}{n}\right)$; when $K>\sqrt{2}\sigma$,
  $\mathbb{E}\left[W_2^2(\PP_n*\mathcal{N}(0, \sigma^2I_d), \PP*\mathcal{N}(0, \sigma^2I_d))\right]$
  converges with rate $\omega\left(\frac{1}{n}\right)$. Here is an obvious gap between $K <
  \sigma/2$ and $K > \sqrt{2}\sigma$, and our results close this gap. Moreover,
  ~\cite{goldfeld2020convergence} shows that $\mathbb{E}\left[W_2(\PP_n*\mathcal{N}(0, \sigma^2I_d), \PP*\mathcal{N}(0, \sigma^2I_d))\right]$ converges with rate
  $\mathcal{O}\left(\frac{1}{n^{1/4}}\right)$ for any $K$ and $\sigma > 0$. 
  
  \begin{proof}[Proof Idea.] Let us introduce the $\chi^2$-mutual information for a pair of random
  variables $S,Y$ as 
  $$ I_{\chi^2}(S; Y) \eqdef \chi^2(P_{S,Y}\| P_S \otimes P_Y)\,,$$ 
  where $\chi^2(P\|Q) = \int \left(\frac{dP}{dQ}\right)^2dQ - 1$.
   
  We will consider the case where $S\sim \PP$, $Y = S + Z$ with $Z\sim \mathcal{N}(0, \sigma^2)$ independent to $S$. 
According to~\cite{goldfeld2020convergence}, the convergence rate of smoothed empirical measure under $W_2,$ KL-divergence and the $\chi^2$-divergence is closely related to $I_{\chi^2}(S; Y)$:\\
(\textbf{Proposition 6} in ~\cite{goldfeld2020convergence}) If $\PP$ is $K$-subgaussian with $K < \sigma$ and $I_{\chi^2}(S; Y) < \infty$, then 
$$\mathbb{E}\left[W_2^2(\PP_n*\mathcal{N}(0, \sigma^2I_d), \PP*\mathcal{N}(0, \sigma^2I_d))\right] = \mathcal{O}\left(\frac{1}{n}\right).$$
(\textbf{Corollary 2} in ~\cite{goldfeld2020convergence}) If $I_{\chi^2}(S; Y) = \infty$, then for any $\tau < \sigma$, 
$$\mathbb{E}\left[W_2^2(\PP_n*\mathcal{N}(0, \tau^2I_d), \PP*\mathcal{N}(0, \tau^2I_d))\right] = \omega\left(\frac{1}{n}\right).$$

Hence our results follow from the following main technical propositions. 
  
\begin{proposition}\label{prop:upper-bound}
	 When $K < \sigma$, for any $K$-subgaussian $d$-dimensional distribution $\mathbb{P}$, we have $I_{\chi^{2}}(S; Y) < \infty$, where $S\sim \mathbb{P}, Z\sim \mathcal{N}(0, \sigma^{2}I_{d}), S\indep Z$ and $Y = S + Z$.
\end{proposition}
\begin{proposition}\label{prop:lower-bound}
    When $K > \sigma$, there exists some 1-dimensional $K$-subgaussian distribution $\mathbb{P}$ such that $I_{\chi^{2}}(S; Y) = \infty$ for $S\sim \mathbb{P}, Z\sim \mathcal{N}(0, \sigma^{2}), S\indep Z$ and $Y = S + Z$.   
\end{proposition}
We will prove these two propositions in the following two sections separately.

We note that results from~\cite{goldfeld2020convergence} and 
Proposition~\ref{prop:upper-bound} also imply 
that $\mathbb{E}[D_{KL}(P_n*\mathcal{N}(0,
\sigma^2I_d)\|P*\mathcal{N}(0, \sigma^2I_d))]$ and $\mathbb{E}[\chi^2(P_n*\mathcal{N}(0,
\sigma^2I_d)\|P*\mathcal{N}(0, \sigma^2I_d))]$ both converge with rate
$O\left(\frac{1}{n}\right)$. Our second Proposition \ref{prop:lower-bound} implies
that for the special $\PP$ constructed there we have 
\begin{align*}\mathbb{E}[D_{KL}(\PP_n*\mathcal{N}(0, \sigma^2I_d)\|\PP*\mathcal{N}(0, \sigma^2I_d)) &=
\omega\left(\frac{1}{n}\right)\\
\mathbb{E}\chi^2(\PP_n*\mathcal{N}(0,
\sigma^2I_d)\|\PP*\mathcal{N}(0, \sigma^2I_d)) &= \infty\,.
\end{align*}
\end{proof}


\par Next, we give a tight characterization on the $W_2$-convergence rate in dimension $d=1$.
\begin{theorem}\label{thm-1D} [Improved bounds for dimension-1] Fix $K > \sigma > 0$ and let
	\begin{equation}\label{eq:def-alpha} \alpha = {(\sigma^{2} + K^{2})^{2} \over 4(\sigma^{4} + K^{4}) }\,. \end{equation}
	\begin{enumerate}
		\item (Lower Bound) There exist a $K$-subgaussian distribution $\PP$ over
		$\mreals$ and $0< \delta_n \to 0$ such that 
			$$ \limsup_{n\to \infty} n^{\alpha + \delta_n}
			\mathbb{E}\left[W_{2}(\mathbb{P}_{n} * \mathcal{N}(0, \sigma^{2}I_{d}),
			\mathbb{P} * \mathcal{N}(0, \sigma^{2}I_{d}))\right] > 0 \,.$$
   		\item (Upper Bound) There exists a sequence $0< \delta_n \to 0$ such that for any 1-dimensional $K$-subgaussian $\PP$ over $\mreals$ 
		and $n\ge 2$, we have
    \begin{equation}\label{eq-thm2-1}\mathbb{E}\left[W_{2}^{2}(\mathbb{P}*\mathcal{N}(0,
    \sigma^{2}), \mathbb{P}_{n}*\mathcal{N}(0, \sigma^{2}))\right] \le n^{-2\alpha + \delta_n}
    \end{equation}
    \end{enumerate}
\end{theorem}
\begin{remark}
    Our proof gives $\delta_n = \frac{1}{\sqrt[3]{\log\log n}}$. With more delicate analysis we believe that our proof could give $\delta_n = \frac{1}{\sqrt[3]{\log n}}$ in the upper bound.
\end{remark}
\begin{remark}
	According to the Cauchy-Schwarz inequality, we have
	$$\mathbb{E}\left[W_{2}(\mathbb{P}*\mathcal{N}(0, \sigma^{2}), \mathbb{P}_{n}*\mathcal{N}(0, \sigma^{2}))\right]\le \sqrt{\mathbb{E}\left[W_{2}^{2}(\mathbb{P}*\mathcal{N}(0, \sigma^{2}), \mathbb{P}_{n}*\mathcal{N}(0, \sigma^{2}))\right]}.$$
	Therefore, the lower bound part in Theorem \ref{thm-1D} indicates that for any $K$ and $\epsilon > 0$, there exists some $K$-subgaussian distribution $\PP$ and $\sigma > 0$ such that 
	\begin{equation}\label{eq-thm2-2}\mathop{\limsup}_{n\to\infty}n^{2\alpha + 2\delta_n}\mathbb{E}\left[W_{2}^{2}(\mathbb{P}_{n} * \mathcal{N}(0, \sigma^{2}I_{d}), \mathbb{P} * \mathcal{N}(0, \sigma^{2}I_{d}))\right] > 0.\end{equation}
	and upper bound part in Theorem \ref{thm-1D} indicates that
	$$\mathbb{E}\left[W_{2}(\mathbb{P}*\mathcal{N}(0, \sigma^{2}), \mathbb{P}_{n}*\mathcal{N}(0, \sigma^{2}))\right]\le n^{-\alpha + \delta_n/2}.$$
\end{remark}

Finally we provide an upper bound on the convergence of smoothed empirical measures under KL divergence:
\begin{theorem}\label{thm-KL}
	Suppose $\mathbb{P}$ is a $d$-dimensional $K$-subgaussian distribution, then for any $\sigma > 0$, we have 
	$$\mathbb{E}\left[D_{KL}\left(\mathbb{P}_{n} * \mathcal{N}(0, \sigma^{2}I_{d})\big\|\mathbb{P} * \mathcal{N}(0, \sigma^{2}I_{d}\right))\right] = \mathcal{O}\left(\frac{(\log n)^{d}}{n}\right).$$
\end{theorem}

\begin{remark}
	From Proposition \ref{prop:upper-bound} and \ref{prop:lower-bound} and also results from
	\cite{goldfeld2020convergence}, we know that when $\sigma > K$, the convergence rate
	is $\mathcal{O}\left(\frac{1}{n}\right)$. From the above theorem, we know that when
	$\sigma\le K$, there exists a $K$-subgaussian distribution $\PP$ such that the convergence rate is $\omega\left(\frac{1}{n}\right)$, and for any $K$-subgaussian distribution, the convergence rate is always
	$\mathcal{O}\left(\frac{(\log n)^{d}}{n}\right)$. Hence at $K=\sigma$, the KL divergence also
	experiences a change to a non-parametric convergence rate, although with only a
	poly-logarithmic slow-down. As we discussed in Corollary~\ref{cor:lsi} this precludes a general, finite logarithmic-Sobolev constant for a Gaussian mixture $\PP*\calN(0,\sigma^2)$ when $\sigma <
	K$.
\end{remark}


\subsection{Organization of this Paper}\label{sec: organization}
\par In Section 2 we will present the proof of Proposition \ref{prop:upper-bound}. In Section 3 we will present the proof of Proposition \ref{prop:lower-bound}. The proof of the lower bound part and the upper bound part of Theorem \ref{thm-1D} will be presented in Section \ref{thm-lower-bound-part} and \ref{thm-upper-bound-part}. Finally in Section \ref{thm-KL-part}, we will present the proof of Theorem \ref{thm-KL}. Finally in the appendix Section \ref{apx:t2}, we prove that a uniform log-Sobolev inequality constant does not exists for distributions we constructed in \ref{sec-bernoulli}.

\subsection{Notations}\label{sec: notations}
\par Throughout this paper, we use $*$ to denote convolutions of two random variables, $i.e.$ for $X\sim\mathbb{P}, Y\sim\mathbb{Q}, X\indep Y$, we have $X + Y\sim \mathbb{P} * \mathbb{Q}$; we use $\otimes$ to denote the product of two random variables's laws, 
$i.e.$ for $X\sim\mathbb{P}, Y\sim\mathbb{Q}, X\indep Y$, we have $(X, Y)\sim \mathbb{P}\otimes\mathbb{Q}$; we use $\circ$ to denote the composition between a Markov kernel $P_{Y|X}$ and a distribution $P_X$, $e.g.$ for $Y$ generated according to $P_{Y|X}$ with $X$'s prior distribution $P_X$, then $Y\sim P_{Y|X}\circ P_X$, cf.~\cite[Section 2.4]{itbook}.

\par Furthermore, we use $\mathbf{P}(E)$ to denote the probability of event $E$, $\mathbb{E}_{\PP}[\cdot]$ to denote the quantity with respect to distribution $\PP$, and we use $\EE[\cdot]$ to denote the expectation after taking expectation of all random variables inside the bracket. 
For sequences indexed by $n\to \infty$ we denote $A_n = \mathcal{O}(B_n), A_n = \Omega(B_n)$  to denote that $A_n\le CB_n$ and $A_n\ge C B_n$ for some positive constant $C$ possibly dependent on $K, \sigma$ and independent of $n$, respectively. We also, equivalently, use $A_n\lesssim B_n$ and $ A_n\gtrsim B_n$, to denote the $A_n= O(B_n)$ and $A_n = \Omega(B_n)$, respectively. We write $A_n = \Theta(B_n)$ whenever $A_n = \mathcal{O}(B_n)$ and $A_n = \Omega(B_n)$ simultaneously. We use $A_n = \tilde{\mathcal{O}}(B_n)$ to denote that $A_n\le CB_n\cdot \log^ln$ for some positive constant $C, l$. We further use $\|\cdot \|$ to denote Euclidean norm, and use $I_d$ to denote the $d\times d$ identity matrix.
\par We will use $\varphi_{\sigma^{2}}(\cdot)$ to denote the density of $d$-dimensional multivariate normal distribution $\mathcal{N}(0, \sigma^{2}I_{d})$. And for 1-dimensional distributions we use $\Phi_{\sigma}$ to denote the CDF of $\mathcal{N}(0, \sigma^{2})$.

\section{Proof of Proposition ~\ref{prop:upper-bound}}
In this section, we provide a proof of Proposition \ref{prop:upper-bound}. The proof idea is to notice that we can write $I_{\chi^{2}}(S; Y)$ as $\mathbb{E}\left[\chi^{2}\left(\mathcal{N}(S, \sigma^{2}I_{d})\|\mathbb{E}\mathcal{N}(S, \sigma^{2}I_{d})\right)\right]$, then we decompose $\mathbb{R}^d$ into several cubes $c_i$ with finite diameters, and we prove an upper bound of $\mathbb{E}\left[\mathbf{1}_{S\in c_i}\chi^{2}\left(\mathcal{N}(S, \sigma^{2}I_{d})\|\mathbb{E}\mathcal{N}(S, \sigma^{2}I_{d})\right)\right]$ for each $i$.

\begin{proof}
\par We suppose that the distribution $\mathbb{P}$ is $d$-dimensional. Then with $S\sim\mathbb{P}, Z\sim\mathcal{N}(0,
\sigma^{2}I_{d}), S\indep Z$ and $Y = S + Z$, noticing that $\chi^2(P_{SY}\|P_S\otimes P_Y) = \mathbb{E}_S[\chi^2(P_{Y|S}\|P_Y)]$, we have 
\begin{equation}\label{eq1}
		I_{\chi^{2}}(S; Y) = \mathbb{E}\left[\chi^{2}\left(\mathcal{N}(S, \sigma^{2}I_{d})\|\mathbb{E}\mathcal{N}(S, \sigma^{2}I_{d})\right)\right]\le \mathbb{E}\left[\int_{\mathbb{R}^{d}}\frac{\varphi_{\sigma^2}(\my - S)^{2}}{\mathbb{E}\varphi_{\sigma^2}(\my - S)}d\my\right]\lesssim \EE\left[\int_{\mathbb{R}^{d}}\frac{\exp\left(-\|\my -
		S\|_{2}^{2}/\sigma^{2}\right)}{\rho(\my)}d\my\right] 
\end{equation}
where $S\sim\mathbb{P}$ and $\rho(\my) = {\mathbb{E}\exp\left(-\|\my - S\|_{2}^{2}/(2\sigma^{2})\right)}\,$.


Let us decompose $\mreals^d=\bigcup_i c_i$ as a union of cubes of diameter $2$. 
Since $K < \sigma$, we have $\frac{K}{\sigma} < 1$. Hence we can choose small $\delta, \delta' > 0$ such that
\begin{equation*}
	\sqrt{\frac{1}{(1 + \delta)^{2}(1 + \delta')}} > \frac{K}{\sigma}.
\end{equation*}
Notice that, due to the $K$-subgaussianity of $S$, we have~\cite[p. 26]{boucheron2013concentration}
\begin{equation}\label{eq:brt_0}
	\EE[\exp\left({(1+\delta')(1+\delta)^2\over 2\sigma^2} \|S\|^2\right)] < \infty 
\end{equation}

We will use the following lower bounds on $\rho(\my)$:\footnote{Notation $\gtrsim$ and $\lesssim$ in this proof
denote inequalities up to constants that may depend on $K,\sigma,d$ and distribution $P$.}
\begin{align} 
  \rho(\my) &\gtrsim \exp\left(-{1+\delta'\over 2\sigma^2} \|\my\|^2\right)\,, \label{eq:brt_1}\\
   \rho(\my) &\ge \exp\left(-\frac{6}{\sigma^2}\right)\mathbf{P}(S\in c_i) \exp\left(-{3\over 4\sigma^2} \|\my-s\|^2\right) \qquad \forall s \in
   c_i\,.\label{eq:brt_2}
\end{align}
Assuming these inequalities, the proof proceeds as follows. Fix an arbitrary $s\in \mreals^d$ and
notice that from~\eqref{eq:brt_1} whenever $\|\my\|\le (1+\delta)\|s\|$ we have
$$ \rho(\my) \gtrsim \exp\left(-{(1+\delta')(1+\delta)^2\over 2\sigma^2} \|s\|^2\right)$$
which implies that  
\begin{equation}\label{eq:brt_4}
	\int_{\|\my\|\le (1+\delta)\|s\|}  \frac{\exp\left(-\|\my -
		s\|_{2}^{2}/\sigma^{2}\right)}{\rho(\my)}d\my \lesssim
		\exp\left({(1+\delta')(1+\delta)^2\over 2\sigma^2} \|s\|^2\right)\,,
\end{equation}		
since the numerator integrates over $\mreals^d$ to $(\pi \sigma^2)^{d/2}$. On the other hand,
from~\eqref{eq:brt_2} if $s\in c_i$ then 
\begin{equation}\label{eq:brt_3}
	\frac{\exp\left(-\|\my -
		s\|_{2}^{2}/\sigma^{2}\right)}{\rho(\my)} \lesssim \mathbf{P}(S\in c_i)^{-1} \exp\left(
		-{\|\my -s \|^2\over 4\sigma^2}\right)\,.
\end{equation}		
Note also that when $\|\my\| \ge (1+\delta)\|s\|$ we have $\|\my -s \|\ge \delta \|s\|$.
Thus, integrating the right-hand side of~\eqref{eq:brt_3} over $\{\my: \|\my-s\| \ge
\delta\|s\|\}$ we obtain
$$ \mathbf{P}(S\in c_i)^{-1} \int_{\|\my-s\| \ge \delta\|s\|} \exp\left(
		-{\|\my -s \|^2\over 4\sigma^2}\right) \lesssim \mathbf{P}(S\in c_i)^{-1} \mathbf{P}\left[U_d >
		{\delta^2 \|s\|^2 \over \sqrt{2\sigma^2}}\right]\,,$$
where $U_d$ denotes a $\chi^2(d)$ random variable with $d$ degrees of freedom. Using Chernoff
inequality $\mathbf{P}(U_d > r) \le 2^{d\over 2} e^{-r/4}$ we obtain 
$$ \max_{s\in c_i} \mathbf{P}\left[U_d >
		{\delta^2 \|s\|^2 \over \sqrt{2\sigma^2}}\right] \lesssim \exp(-C \|x_i\|^2)\,, $$
where $x_i$ is the center of the cube $c_i$ and $C$ is some constant.

Thus, together with~\eqref{eq:brt_4} we obtain that for any $s\in c_i$:
$$ \chi^2(\mathcal{N}(s, \sigma^2 I_d) \| P_Y) \lesssim \int_{\mreals^d}  \frac{\exp\left(-\|\my -
		s\|_{2}^{2}/\sigma^{2}\right)}{\rho(\my)}d\my \lesssim \mathbf{P}(S\in c_i)^{-1} \exp(-C
		\|x_i\|^2) + \exp\left({(1+\delta')(1+\delta)^2\over 2\sigma^2} \|s\|^2\right)\,.$$

Taking expectation of the latter over $S$, the second term is finite because of~\eqref{eq:brt_0},
while the first one is bounded because the number of cubes with $\|x_i\| \le r$ is $O(r^d)$. This
completes the proof of finiteness of~\eqref{eq1}, assuming~\eqref{eq:brt_1} and~\eqref{eq:brt_2}.
We now establish these.

To show~\eqref{eq:brt_1}  set $t$ to be any value such that $\mathbf{P}(\|S\|<t) \ge {1\over 2}$ and notice that
\begin{equation}\label{eq:brt_5}
	\rho(\my) \gtrsim \EE[\exp\left(-{\|\my - S\|^{2}\over 2\sigma^{2}}\right) | \|S\|<t ] \,.
\end{equation}
Next, notice that for any $t$ and $\delta'>0$ we can find some constant $C'$
such that 
\begin{equation}\label{eq:brt_6}
	(a+t)^2 \le (1+\delta')a^2 + C',\qquad \forall a\in\mreals\,.
\end{equation}
Thus for any $\|s\|<t$ we have
$$ \exp\left(-{\|\my - s\|^{2}\over2\sigma^2}\right) \ge
\exp\left(-{(\|\my\| + \|s\|)^{2}\over2\sigma^{2}}\right) \gtrsim \exp\left(-{\|\my\|^2
(1+\delta')\over2\sigma^{2}}\right)\,. $$
Using this estimate in~\eqref{eq:brt_5} recovers~\eqref{eq:brt_1}. 

To show~\eqref{eq:brt_2} we start similarly:
\begin{equation}\label{eq:brt_7}
	\rho(\my) \ge \mathbf{P}(S\in c_i) \EE[\exp\left(-{\|\my - S\|^{2}\over 2\sigma^{2}}\right) | S\in c_i ].
\end{equation}
Now fix any (non-random) $s\in c_i$ and notice that under the conditioning above we have
$\|S-s\|\le 2$ because the cube $c_i$ has diameter $2$. Then from triange inequality and
~\eqref{eq:brt_6} with $\delta'=1/2$ we obtain
$$ \|\my-S\|^2 \le (\|\my-s\| + 2)^2 \le {3\over 2} \|\my-s\|^2 + 12\,.$$
Using this bound in~\eqref{eq:brt_7} yields~\eqref{eq:brt_2}.
\end{proof}

In some of our applications in the sequel, we also require that the R\'{e}nyi mutual information $I_\lambda(S; Y)$ (defined in Definition \ref{def-renyi}) is finite for all $1<\lambda < 2$, which we demonstrate using a similar approach in Lemma \ref{lem-I-lambda}.

\section{Proof of Proposition~\ref{prop:lower-bound}}\label{sec-lower-bound}
\par In this section, we will present a proof of Proposition \ref{prop:lower-bound}. The main idea of this proof is to construct a hard example $\mathbb{P} = \sum_{k=0}^{\infty} p_{k}\delta_{r_{k}}$ with subgaussian tails, where $r_i$ and $r_j$ are far away from each other so that $\delta_{r_j} * \mathcal{N}(0, \sigma^2)$ with $j\neq i$ has very little effects on the density of $\PP * \mathcal{N}(0, \sigma^2)$ near $r_i$. Therefore we can uniformly lower bound $\mathbb{E}\left[\mathbf{1}_{S = r_i}\chi^{2}\left(\mathcal{N}(S, \sigma^{2}I_{d})\|\mathbb{E}\mathcal{N}(S, \sigma^{2}I_{d})\right)\right]$ for each $i$, and if we sum up over all $i$ we can prove the infiniteness of $I_{\chi^2}(S; Y)$. 
\begin{proof}
Without loss of generality we assume $\sigma = 1$, and we only need to prove the proposition for $K > 1$. (Otherwise we consider $S' = S/\sigma, Z' = Z/\sigma$ and $Y' = Y/\sigma$, and we will have $S'$ is a $K/\sigma$-subgaussian distribution, $Z'\sim \mathcal{N}(0, 1)$ and $I_{\chi^{2}}(S; Y) = I_{\chi^{2}}(S'; Y')$.)

We construct a 1-dimensional distribution $\mathbb{P}$ similarly to~\cite{goldfeld2020convergence} as $\mathbb{P} = \sum_{k=0}^{\infty} p_{k}\delta_{r_{k}}$, where we choose $r_{0} = 0$, $p_{0} = 1 - \sum_{k=1}^{\infty}p_{k}$ and for some positive constant $c_{1}$ to be determined we choose
\begin{equation}\label{eq-prob}
	p_{k} = c_{1}\exp\left(-\frac{r_{k}^{2}}{2K^{2}}\right), \quad k\ge 1.
\end{equation}
Here we let $r_{i}$ be a geometric sequence: $r_{k} = c^{k-1}$ for any $k\ge 1$, with $c > 2$ is a constant to be specified later. We restrict that $c_{1}$ only depends on $K$ and $c_{1}\cdot \sum_{k=1}^{\infty}\exp(- r_{k}^{2}/(2K^{2})) < 1.$ Then we will have $p_{0} = 1 - \sum_{k=1}^{\infty}p_{k} > 0$ making $\PP$ a well-defined
distribution.

\par 
To show that our distribution is $K$-subgaussian, we need to verify the
definition~\eqref{eq:subgauss}. Here we only show a weaker claim that 
    \begin{equation}\label{eq:subgauss_weak}
    	\mathbb{E}[\exp\left(\alpha S\right)] \le 2\exp\left(\frac{K^2\alpha^2}{2}\right).
	\end{equation}    
   The full proof, found in~\cite[Appendix A]{block2022rate}, is removed to save space.

To that end, we notice that
    $$\mathbb{E}[\exp\left(\alpha S\right)] = p_{0} + c_{1}\sum_{k=1}^{\infty}\exp\left(-\frac{1}{2K^{2}}\left(r_{k} - \alpha K^{2}\right)^{2}\right)\exp\left(\frac{K^{2}\alpha^{2}}{2}\right).$$
    We suppose $k_0$ to be the smallest $k$ such that $r_k - \alpha K^{2}$ is positive. Since
    $r_{k+1} - r_{k}\ge 1$ for every $k$, we have for $k\ge k_{0}$, $r_{k} - \alpha K^{2}\ge k -
    k_{0} + r_{k_0} - \alpha K^{2}\ge k - k_{0}$, and for $k < k_0$, $r_{k} - \alpha K^2\le r_{k_0 -
    1} - \alpha K - (k_{0} - 1 - k)\le -(k_{0} - 1 - k)$ since $r_{k_0 - 1}-\alpha K^2\le 0$. Hence, we have $\sum_{k=1}^{\infty} e^{-\frac{1}{2K^{2}}\left(r_{k} - \alpha K^{2}\right)^{2}} \le 2 \sum_{k=0}^{\infty}e^{-\frac{k}{2K^2}}= \frac{2}{1 - \exp\left(-\frac{1}{2K^2}\right)}.$
    Therefore, if we choose $c_1 = (1 - \exp(-1/(2K^2))/2$, and notice that $p_0\le 1\le \exp(K^2\alpha^2/2)$, we would have
    $$\mathbb{E}[\exp\left(\alpha S\right)] \le \exp\left(\frac{K^2\alpha^2}{2}\right) + \exp\left(\frac{K^2\alpha^2}{2}\right) = 2\exp\left(\frac{K^2\alpha^2}{2}\right).$$

For now we proceed assuming that $c_1$ is already chosen such that for any $c>2$, we have that $\mathbb{P}$
is a $K$-subgaussian distribution. Then, our goal is to specify a value of constant $c$ such that $I_{\chi^{2}}(S; Y) = \infty$. 
From the definition of $I_{\chi^{2}}$, we have the following chain: for any fixed $\delta>0$,
\begin{align*}
	I_{\chi^{2}}(S; Y) + 1 &\eqdef \int_{\mathbb{R}}\frac{\mathbb{E}\varphi_{1}^{2}(y - S)}{\mathbb{E}\varphi_{1}(y - S)}dy \gtrsim \int_{\mathbb{R}}\frac{\sum_{k=0}^{\infty}p_{k}\varphi_{\frac{1}{\sqrt{2}}}(y - r_{k})}{\sum_{k=1}^{\infty}p_{k}\varphi_{1}(y - r_{k})}dy = \sum_{k=0}^{\infty}\int_{\mathbb{R}}\frac{p_{k}\varphi_{\frac{1}{\sqrt{2}}}(y - r_{k})}{\sum_{k=1}^{\infty}p_{k}\varphi_{1}(y - r_{k})}dy\\
	&\gtrsim \sum_{k=0}^{\infty}\int_{r_k-\delta}^{r_k+\delta}\frac{\varphi_{\frac{1}{\sqrt{2}}}(y - r_{k})}{\varphi_{1}(y - r_{k})}\cdot\left(1 + \sum_{j\neq k}\frac{p_{j}}{p_{k}}\frac{\varphi_{1}(y - r_{j})}{\varphi_{1}(y - r_{k})}\right)^{-1}dy\,. \numberthis \label{eq:brt_8}
\end{align*}
where in the last inequality we use the fact that the integrand is always nonnegative. Next, we will upper bound the term $1 + \sum_{j\neq k}\frac{p_{j}}{p_{k}}\frac{\varphi_{1}(y - r_{j})}{\varphi_{1}(y - r_{k})}$ in the above inequality.

First, consider $j=0$ and $|y-r_k|\le
\delta$. From the definition of $p_k$ and $p_0\le 1$ we have
\begin{equation}\label{eq: decomposoition j0}
	\log \frac{c_1p_0\varphi_{1}(y)}{p_{k}\varphi_{1}(y - r_{k})} \le -\frac{y^{2}}{2} + \frac{r_{k}^{2}}{2K^{2}} + \frac{(y - r_{k})^{2}}{2} \le -\frac{(r_{k} - \delta)^{2}}{2} + \frac{r_{k}^{2}}{2K^{2}} + \frac{\delta^{2}}{2}\,.
\end{equation}
When $K>1$ the latter is upper-bounded by $\frac{\delta^{2}}{2K^{2}(1 - K^{2})} +
\frac{\delta^{2}}{2}$ uniformly for all $k$ and $y$. Thus, we get that term 
$ \frac{p_{0}}{p_{k}}\frac{\varphi_{1}(y - r_{0})}{\varphi_{1}(y - r_{k})}  \le C$ for some
$C=C(K)$.

For $j\ge 1$ and $|y - r_{k}|\le \delta$, we have by bounding $y(r_j-r_k) \le -r_k^2 + r_kr_j + \delta|r_k-r_j|\le -r_k^2 + r_kr_j + \delta r_k + \delta r_j$ the following chain
\begin{align*}
	&\quad \log \frac{p_{j}}{p_{k}}\frac{\varphi_{1}(y - r_{j})}{\varphi_{1}(y - r_{k})} = \left(\frac{1}{2K^{2}} + \frac{1}{2}\right)(r_{k}^{2} - r_{j}^{2}) - y(r_{k} - r_{j})\\
	&\le \left(\frac{1}{2K^{2}} + \frac{1}{2}\right)(r_{k}^{2} - r_{j}^{2}) -r_k^2 + r_kr_j + \delta r_k + \delta r_j = A_1 + A_2 + A_3-\frac{r_j^2}{4},
\end{align*}
where we denoted $A_1\eqdef \frac{l}{2}r_{k}^{2} - \frac{1}{2K^{2}} r_{j}^{2} + r_{k}r_{j}, A_2\eqdef {\ell\over 2}r_k^2 +\delta r_k$ and $A_3\eqdef -{1\over 4}r_j^2 + \delta r_j$ with $\ell \eqdef \frac{1}{2K^{2}} - \frac{1}{2}$.

Note that $K>1$ and, thus, $\ell < 0$.  We show that by choosing $c$ and $\delta$ it is possible to make sure $A_1, A_2, A_3\le 0$ for all $k,j$. First, notice that because $r_k\ge 1$ or $r_k=0$ by setting $\delta = \min\left(-{\ell\over 2}, {1\over 4}\right)$ we have $A_2, A_3\le 0$. Second, we have $A_1 = r_j^2 f(r_k/r_j)$ where $f(y) = {\ell \over 2} y^2 + y -{1\over 2K^2}$. Since $f(0)<0$ and $f(+\infty)=-\infty$ (recall $\ell < 0$) we must have that for some sufficiently large $c>0$ we have $f(y) < 0$ if $y \le 1/c$ or $y \ge c$. For convenience we take this $c>2$ as well.  Since $r_k/r_j$ is always either $\le 1/c$ or $\ge c$ we conclude $A_1 \le 0$. 

Continuing, we obtained that with our choice of $c$, for $j\neq k, j\ge 1$ and $|y - r_{k}|\le \delta$ we have
\begin{align*}
	\log \frac{p_{j}}{p_{k}}\frac{\varphi_{1}(y - r_{j})}{\varphi_{1}(y - r_{k})} \le A_1 + A_2 + A_3 - {r_j^2\over 4}\le  - \frac{r_{j}^{2}}{4}  \le   -\frac{c^{j}}{4}\lesssim -\frac{j}{2},
\end{align*}
Combining this inequality with \eqref{eq: decomposoition j0}, we obtain that 
\begin{equation*}
	1 + \sum_{j\neq k}\frac{p_{j}}{p_{k}}\frac{\varphi_{1}(y - r_{j})}{\varphi_{1}(y - r_{k})}\le 1 + C + \sum_{j = 1, j\neq k}^{\infty}\exp(-j/2)\le 1 + C + \frac{1}{1 - \exp(-1/2)}\triangleq C'.
\end{equation*}

Bringing this inequality back to \eqref{eq:brt_8}, we obtain that 
\begin{align*}
	\sum_{k=0}^{\infty}\int_{r_{k} - \delta}^{r_{k} + \delta}\frac{\varphi_{\frac{1}{\sqrt{2}}}(y - r_{k})}{\varphi_{1}(y - r_{k})}\cdot\left(1 + \sum_{j\neq k}\frac{p_{j}}{p_{k}}\frac{\varphi_{1}(y - r_{j})}{\varphi_{1}(y - r_{k})}\right)^{-1}dy\ge \left(\int_{ - \delta}^{\delta}\frac{\varphi_{\frac{1}{\sqrt{2}}}(y)}{\varphi_{1}(y)}dy\right)\cdot \sum_{k=0}^{\infty}\frac{1}{C'} =  \infty
\end{align*}
And we have proved that $I_{\chi^{2}}(S; Y) = \infty$.
\end{proof}

\section{Proof of the Lower Bound in Theorem \ref{thm-1D}}\label{thm-lower-bound-part}

The main idea of the lower bound is a simple observation encoded in the next
proposition (proved in
Appendix~\ref{sec-bernoulli}). It turns out that already simple binary Gaussian mixtures
demonstrate why for $K>\sigma$ the rate of smoothed Wasserstein convergence should slow down. 

\begin{proposition}\label{prop-bernoulli}
	Fix $K>\sigma > 0$. For every $h > 0$ we consider a Bernoulli
	distribution $\mathbb{P}_{h} = (1-p_h)\delta_0 + p_h\delta_h$, with $p = e^{-h^2/(2K^2)}$.
	Then for any $\epsilon > 0$ we have
    $$\sup_h\mathbb{E}\left[W_{2}(\mathbb{P}_{h} * \mathcal{N}(0, \sigma^{2}), \mathbb{P}_{h, n} * \mathcal{N}(0, \sigma^{2}))\right] = \Omega\left(n^{-\alpha-\epsilon}\right),$$
 	where $\mathbb{P}_{h, n}$ is the empirical measure constructed from $n$ $i.i.d.$ samples from $\mathbb{P}_{h}$, and $\alpha$ is defined in \eqref{eq:def-alpha}.
\end{proposition}

Note that the rate in the lower bound coincides with the rate in the lower bound in Theorem
\ref{thm-1D}. So to prove Theorem \ref{thm-1D} we only need to show that there exists a single
$\mathbb{P}$ which works for all $n$ (in the above the choice of $h$ implicitly depends on $n$).
To achieve that goal, our proof below takes a convex combination of many Bernoulli distributions
to produce a single $K$-subgaussian distribution that achieves the rate.

\par We construct the following discrete distribution
\begin{equation}
	\mathbb{P} = \sum_{k=0}^{\infty} p_{k}\delta_{r_{k}},\quad p_{k}\ge 0, \quad\sum_{k=1}^{\infty}p_{k} = 1,
\end{equation}
where we choose $r_{k} = c_1c_2c_3\cdots c_{k-1}$ for $k\ge 1$ for some positive constant $3\le c_1\le c_2\le c_3\le \cdots$ to be determined later, and
\begin{equation}\label{eq: def-p-k}
p_{k} = \frac{C}{\sqrt{2\pi}K}\exp\left(- \frac{r_{k}^{2}}{2K^{2}}\right), \quad k\ge 1,\quad p_0 = 1 - \sum_{k=1}^\infty p_k
\end{equation}
where $C$ is a small enough constant such that $\sum_{k=1}^\infty p_k\le 1$. We claim that
$\mathbb{P}$ is $K$-subguassian (in the sense of~\eqref{eq:subgauss}), but we omit the details of
this verification since a very similar distribution was already analyzed in Proposition
\ref{prop:lower-bound}.

\par We let $\kappa = \frac{\sigma^{2}}{K^{2}}\in (0, 1)$, and 
\begin{equation}\label{eq: def-t-k}
	t_k = \frac{1}{2}(c_k+1)\left(1 + \kappa\right)\ge \frac{1}{2}(c_k + 1)\ge 2.
\end{equation}
First we will provide the following proposition, which upper and lower bounds the probability of $\mathbb{P} * \mathcal{N}(0, \sigma^{2})$ near $t_kr_{k}$. The proof of this proposition is deferred to the Appendix \ref{sec: appendix-proofs}.

\begin{proposition} \label{prop: lb-lb-ub} There exist $C_l,C_u > 0$ such that for all $k\ge 1$ we have
\begin{align*}
	2\sigma^2 \log \mathbf{P}(X\in [t_kr_{k} + 1, t_kr_{k} + 2]) &\ge -\left(t_k^{2} - \kappa c_k - c_k\right) (r_{k} + 2)^{2} - C_l \numberthis \label{eq: lb-lb}\\
	2\sigma^2 \log \mathbf{P}(X\in [t_kr_{k} , t_kr_{k} + 2]) &\le -\left(t_k^{2} - \kappa
	c_k - c_k\right) (r_{k} - 2)^{2} + C_u\,, \numberthis \label{eq: lb-ub}
\end{align*}
where the latter holds whenever $c_k\ge \max\left\{\sqrt{\frac{2}{\kappa}}, \frac{\kappa + 3}{1 -
\kappa}\right\}$.
\end{proposition}

\par We next present the following proposition, indicating that with positive probability the difference of CDFs of $\PP*\mathcal{N}(0, \sigma^{2})$ and $\PP_{n}*\mathcal{N}(0, \sigma^{2})$ is larger than $\frac{1}{2}\sqrt{\frac{p_{k+1}}{n}}$, which we will show is, in turn, larger than $\mathbb{P}*\mathcal{N}(0, \sigma^{2})([t_kr_{k}, t_kr_{k} + 2])$ under some assumptions.
\begin{proposition}\label{prop-lb-difference}
	Suppose $c_k\ge \frac{\kappa + 3}{1 - \kappa}$ for every $k$. We use $F_\sigma$ and $\tilde{F}_{n, \sigma}$ to denote the CDF of $P*\mathcal{N}(0, \sigma^{2})$ and $P_{n}*\mathcal{N}(0, \sigma^{2})$ respectively. Then there exists $k_{0} = k_0(\sigma, K, C) > 0$ and some global constant $C' > 0$ such that $\forall k\ge k_{0}$ and $n$ with $np_{k+1}\ge C'$, with probability at least $\frac{1}{64}$ we have
$$\tilde{F}_{n, \sigma}(t_kr_{k}) - F_\sigma(t_kr_{k})\ge \frac{1}{2}\sqrt{\frac{p_{k+1}}{n}}.$$
\end{proposition}
\begin{proof}
First we can write
\begin{align*}
	F_\sigma(t_kr_{k}) = \sum_{j=0}^{\infty}p_{j}\Phi_{\sigma}(t_kr_{k} - r_{j}),\quad \text{and}\quad \tilde{F}_{n, \sigma}(t_kr_{k}) = \sum_{j=0}^{\infty}\hat{p}_{j}\Phi_{\sigma}(t_kr_{k} - r_{j}),
\end{align*}
where $\Phi_{\sigma}$ is CDF of $\mathcal{N}(0, \sigma^{2})$, and $\hat{p}_{j}$ is the empirical estimation of $p_{j}$ with these $n$ samples. Then we have
\begin{align*}
\tilde{F}_{n, \sigma}(t_kr_{k}) - F_\sigma(t_kr_{k}) & = \sum_{j=0}^{\infty}(\hat{p}_{j} - p_{j})\Phi_{\sigma}(t_kr_{k} - r_{j})\\
&\ge \sum_{j=0}^{k}\hat{p}_{j} - \sum_{j=0}^{k}p_{j} - \sum_{j=0}^{k}|\hat{p}_{j} - p_{j}|(1 - \Phi_{\sigma}(t_kr_{k} - r_{j})) - \sum_{j=k+1}^{\infty}|\hat{p}_{j} - p_{j}|\Phi_{\sigma}(t_kr_{k} - r_{j}).
\end{align*}
From assumption $c_k\ge \frac{\kappa + 3}{1 - \kappa}$ we know that $c_k\ge t_k + 1$. Hence for any $j\ge k+1$ we have $|t_kr_{k} - r_{j}|\ge |(c_k - t_k)r_{k}|\ge r_{k}\ge 1$ and for any $j\le k$ we have $|t_kr_{k} - r_{j}|\ge (t_k - 1)r_{j}\ge r_{j}\ge 1$. According to the upper bound of Gaussian tail function (Proposition 2.1.2 in \cite{vershynin2018high}), we have for $j$ with $t_kr_k - r_j > 0$,
\begin{align*}
	1 - \Phi_{\sigma}(t_kr_{k} - r_{j}) & \le \frac{1}{\sqrt{2\pi}}\cdot \frac{\sigma}{|t_kr_k - r_j|}\exp\left(-\frac{(t_kr_{k} - r_{j})^{2}}{2\sigma^{2}}\right)\le \sigma\exp\left(-\frac{(t_kr_{k} - r_{j})^{2}}{2\sigma^{2}}\right),
\end{align*}
and similarly for $j$ with $t_kr_k - r_j < 0$, $\Phi_{\sigma}(t_kr_{k} - r_{j})\le \sigma\exp\left(-\frac{(t_kr_{k} - r_{j})^{2}}{2\sigma^{2}}\right)$.

\par In the next, we choose $C' = 32768$, and given that $np_{k+1}\ge C'$, we will provide both a lower bound to $\sum_{j=0}^{k}\hat{p}_{j} - \sum_{j=0}^{k}p_{j}$ and also an upper bound to $|\hat{p}_{k+1} - p_{k+1}|$. As for $\sum_{j=0}^{k}\hat{p}_{j} - \sum_{j=0}^{k}p_{j}$, we can write it as
$$\sum_{j=0}^{k}\hat{p}_{j} - \sum_{j=0}^{k}p_{j}  = \frac{1}{n}\left(\sum_{l=1}^{n}U_{l}\right) - \mathbb{E}[U_{1}],$$
 where $U_{l}\sim \text{Bern}(\sum_{j=k+1}^{\infty}p_{j})$ are \emph{i.i.d.} Bernoulli random variables. According to the Berry-Esseen Theorem~\citep{durrett2019probability} 
 we have 
$$\left|\mathbf{P}\left(\frac{1}{\sqrt{n\text{Var}[U_{1}]}}\sum_{l=1}^{n}[U_{l} - \mathbb{E}U_{1}]\ge 1\right) - \mathbf{P}(V\ge 1)\right|\le \frac{\mathbb{E}|U_{1} - \mathbb{E}[U_{1}]|^{3}}{2\sqrt{n}\sqrt{\text{Var}[U_{1}]}^{3}}$$
where $V\sim \mathcal{N}(0, 1)$. It is easy to check that $\sum_{j=k+1}^{\infty}p_{j}\le 2p_{j+1} < 1/2$ for $k\ge 2$. Hence we have
\begin{align*}
	\text{Var}[U_{1}] \ge \frac{1}{2}\sum_{j=k+1}^{\infty}p_{j}\ge \frac{1}{2}p_{k+1}\quad \text{and}\quad \mathbb{E}|U_{1} - \mathbb{E}[U_{1}]|^{3} \le \mathbb{E}|U_{1}|^{3} = \mathbb{E}[U_{1}] = \sum_{j=k+1}^{\infty}p_{j}\le 2p_{k+1}.
\end{align*}
Noticing that for standard Gaussian random variable $V\sim \mathcal{N}(0, 1)$ we have $P(V > 1)\ge 1/8$, we obtain that 
\begin{align*}
	\mathbf{P}\left(\sum_{j=0}^{k}\hat{p}_{j} - \sum_{j=0}^{k}p_{j}\ge \sqrt{\frac{p_{k+1}}{2n}}\right) & = \mathbf{P}\left(\frac{1}{n}\sum_{l=1}^{n}U_{l} - \mathbb{E}[U_{1}]\ge \sqrt{\frac{p_{k+1}}{2n}}\right) \ge \mathbf{P}\left(\frac{1}{\sqrt{n\text{Var}[U_{1}]}}\sum_{l=1}^{n}U_{l} - \mathbb{E}[U_{1}]\ge 1\right)\\
	& \ge \frac{1}{8} - \frac{\mathbb{E}|U_{1} - \mathbb{E}[U_{1}]|^{3}}{2\sqrt{n}\sqrt{\text{Var}[U_{1}]}^{3}}\ge \frac{1}{8} - \frac{2\sqrt{2}}{\sqrt{np_{k+1}}}\ge \frac{1}{16},
\end{align*}
where we use the fact that $np_{k+1}\ge 32768$. As for $|\hat{p}_{k+1} - p_{k+1}|$, if we let $U_{l}'\sim \text{Bern}(p_{k+1}), i.i.d$, again according to Berry-Esseen~\citep{durrett2019probability} Theorem we obtain that
\begin{align*}
	\mathbf{P}\left(|\hat{p}_{k+1} - p_{k+1}| \ge 8\sqrt{\frac{p_{k+1}}{n}}\right) & = \mathbf{P}\left(\frac{1}{n}\sum_{l=1}^{n}U'_{l} - \mathbb{E}[U'_{1}]\ge 8\sqrt{\frac{p_{k+1}}{n}}\right)\le \mathbf{P}\left(\frac{1}{\sqrt{n\text{Var}[U_{1}']}}\sum_{l=1}^{n}U_{l}' - \mathbb{E}[U_{1}']\ge 8\right)\\
	& \le \frac{1}{128} + \frac{\mathbb{E}|U_{1}' - \mathbb{E}[U_{1}']|^{3}}{2\sqrt{n}\sqrt{\text{Var}[U_{1}']}^{3}}\le \frac{\sqrt{2}\mathbb{E}|U_{1}'|^{3}}{\sqrt{n}\sqrt{p_{k+1}}^{3}}\le \frac{1}{128} + \frac{\sqrt{2}}{\sqrt{np_{k+1}}}\le \frac{1}{64}
\end{align*}
after noticing that $\text{Var}[U_{1}] = p_{k+1}(1-p_{k+1})\le p_{k+1}$ and also $P(V > 8)\le 1/128$ for $V\sim \mathcal{N}(0, 1)$, and the last inequality follows from $np_{k+1}\ge 32768$.

\par We further notice that $\mathbb{E}\left[\max_{j\ge 0}|\hat{p}_{j} - p_{j}|^{2}\right]\le \mathbb{E}\left[\sum_{j = 0}^{\infty}|\hat{p}_{j} - p_{j}|^{2}\right] = \sum_{j=0}^{\infty}p_{j}(1-p_{j})/n\le 1/n$. Hence by the Markov inequality we obtained that $\mathbf{P}\left(\max_{j\ge 0}|\hat{p}_{j} - p_{j}|\le 4/\sqrt{n}\right)\ge 15/16$. Therefore, if $n\ge 32768/p_{k+1}$,  with probability at least $1/64$ we have
\begin{align*}
	\tilde{F}_{n, \sigma}(t_kr_{k}) - F_\sigma(t_kr_{k}) & \ge \sqrt{\frac{p_{k+1}}{2n}} - \frac{4\sigma}{\sqrt{n}}\sum_{j=0}^{k}\exp\left(-\frac{(t_kr_{k} - r_{j})^{2}}{2\sigma^{2}}\right) - \frac{4\sigma}{\sqrt{n}}\sum_{j=k+2}^{\infty}\exp\left(-\frac{(t_kr_{k} - r_{j})^{2}}{2\sigma^{2}}\right)\\
	&\quad - 8\sigma\sqrt{\frac{p_{k+1}}{n}}\exp\left(-\frac{(t_kr_{k} - r_{k+1})^{2}}{2\sigma^{2}}\right). \numberthis \label{eq8}
\end{align*}

\par Additionally, we have $\sum_{j=0}^{k}\exp\left(-\frac{(t_kr_{k} - r_{j})^{2}}{2\sigma^{2}}\right)\le k\exp\left(-\frac{(t_k-1)^{2}r_{k}^{2}}{2\sigma^{2}}\right).$ And for any $j\ge k+2$, we have $r_{j} - t_kr_{k}\ge j - (k+2) + r_{k+2} - t_kr_{k}\ge j - (k+2) + (t_k-1)t_{k}$, which indicates that
$$\sum_{j=k+2}^{\infty}\exp\left(-\frac{(t_kr_{k} - r_{j})^{2}}{2\sigma^{2}}\right)\le \left(\sum_{j=k+2}^{\infty}\exp\left(-\frac{j - (k+2)}{2\sigma^{2}}\right)\right)\cdot \exp\left(-\frac{(t_k-1)^{2}r_{k}^{2}}{2\sigma^{2}}\right)\le C_{j}\exp\left(-\frac{(t_k-1)^{2}r_{k}^{2}}{2\sigma^{2}}\right),$$
where $C_{j}$ is a constant only depending on $\sigma$. We also notice that $\frac{(t_kr_{k} - r_{k+1})^{2}}{2\sigma^{2}}\ge \frac{r_{k}^{2}}{2\sigma^{2}}$ using the fact that $c_k\ge t_k + 1$, and that
$$\exp\left(-\frac{(t_k-1)^{2}r_{k}^{2}}{2\sigma^{2}}\right)\le \exp\left(-\frac{c_k^{2}r_{k}^{2}}{4K^{2}} - \frac{c_k^{2}r_{k}^{2}\kappa^{2}}{8\sigma^{2}}\right) = \sqrt{\frac{\sqrt{2\pi}Kp_{k+1}}{C}}\cdot \exp\left(-\frac{c_k^{2}\kappa^{2}r_{k}^{2}}{8\sigma^{2}}\right)$$
using the fact that $2c_k^{2}\kappa + c_k^{2}\kappa^{2}\le c_k^{2}\kappa^{2} + c_k^{2} + \kappa^{2} + 1 - 2c_k - 2\kappa + 2c_k^{2}\kappa = (2t_k - 2)^{2}$. Hence we have
\begin{align*}
	&\quad \frac{4\sigma}{\sqrt{n}}\sum_{j=0}^{k}\exp\left(-\frac{(t_kr_{k} - r_{j})^{2}}{2\sigma^{2}}\right) + \frac{4\sigma}{\sqrt{n}}\sum_{j=k+2}^{\infty}\exp\left(-\frac{(t_kr_{k} - r_{j})^{2}}{2\sigma^{2}}\right) + 8\sigma\sqrt{\frac{p_{k+1}}{n}}\exp\left(-\frac{(t_kr_{k} - r_{k+1})^{2}}{2\sigma^{2}}\right)\\
	& \le 4\sqrt{\frac{p_{k+1}}{n}}\cdot \sigma\left(\frac{C_j\sqrt{\sqrt{2\pi}K}+k}{\sqrt{C}}\exp\left(-\frac{c_k^{2}\kappa^{2}r_{k}^{2}}{8\sigma^{2}}\right) +\exp\left(-\frac{r_{k}^{2}}{2\sigma^{2}}\right)\right).
\end{align*}
Since $r_{k} = c_1c_2\cdots c_{k-1}\ge 3^{k-1}$, there exists some constant $k_{0}$ only depending on $K, \sigma, C$ such that for any $k\ge k_{0}$, we have
$$\sigma\left(\frac{C_j\sqrt{\sqrt{2\pi}K}+k}{\sqrt{C}}\exp\left(-\frac{c_k^{2}\kappa^{2}r_{k}^{2}}{8\sigma^{2}}\right) +\exp\left(-\frac{r_{k}^{2}}{2\sigma^{2}}\right)\right)\le \frac{1}{4\sqrt{2}} - \frac{1}{8}$$
Bringing this result to \eqref{eq8}, we will obtain that for any $k\ge k_{0}$,
$$\tilde{F}_{n, \sigma}(t_kr_{k}) - F_\sigma(t_kr_{k})\ge \frac{1}{2}\sqrt{\frac{p_{k+1}}{n}}$$
holds. This completes the proof of this proposition. 
\end{proof}
With the above propositions, we are now ready to prove the lower bound part of Theorem \ref{thm-1D}.

\begin{proof}[Proof of Lower Bound Part of Theorem \ref{thm-1D}] 
We let $t_k = \frac{1}{2}(1 + \kappa)(1 + c_k)$ and 
\begin{equation}\label{eq7}n_{k} = \left\lfloor\frac{1}{4C_{u}^{2}}\exp\left((t_k^{2} - c_k\kappa - c_k)\cdot \frac{(r_{k} - 2)^{2}}{\sigma^{2}} - \frac{c_k^{2}r_{k}^{2}}{2K^{2}}\right)\right\rfloor,\end{equation}
Then there exists some constant $k_{0}'$ only depending on $k, \sigma$ and $C$ such that for any $k\ge k_{0}'$, we would have $n_kp_{k+1}\ge 32768$. 
Hence according to Proposition \ref{prop-lb-difference} we would have when $k\ge \max\{k_{0}, k_{0}'\}$,
$$\tilde{F}_{n_k, \sigma}(t_kr_{k}) - F_\sigma(t_kr_{k})\ge \frac{1}{2}\sqrt{\frac{p_{k+1}}{n_k}}$$
holds with probability at least $\frac{1}{64}$. Moreover, with our choice of $n_k$, 
it is easy to check that
$$\frac{1}{2}\sqrt{\frac{p_{k+1}}{n_k}}\ge C_{u}\exp\left(-(t_k^{2} - c_k\kappa - c_k)\cdot \frac{(r_{k} - 2)^{2}}{2\sigma^{2}}\right).$$
Hence according to Proposition \ref{prop: lb-lb-ub}, with probability at least $\frac{1}{64}$ we have for $X\sim \mathbb{P}*\mathcal{N}(0, \sigma^{2})$,
$$\tilde{F}_{n_k, \sigma}(t_kr_{k}) - F_\sigma(t_kr_{k})\ge C_{u}\exp\left(-(t_k^{2} - c_k\kappa - c_k)\cdot \frac{(r_{k} - 2)^{2}}{2\sigma^{2}}\right)\ge \mathbf{P}(X\in [t_kr_{k}, t_kr_{k} + 2]).$$
Therefore we have $\tilde{F}_{n_k, \sigma}(t_kr_{k})\ge F_\sigma(t_kr_{k} + 2).$ According to Lemma \ref{lem2} and Proposition \ref{prop: lb-lb-ub}, this indicates that with probability at least $\frac{1}{64}$, 
\begin{align*}
	W_{2}(\mathbb{P} * \mathcal{N}(0, \sigma^2), \mathbb{P}_{n_k} * \mathcal{N}(0, \sigma^2)) & \ge \sqrt{\mathbf{P}(X\in [t_kr_{k} + 1, t_kr_{k} + 2])}\\
	& \ge \sqrt{C_l\exp\left(-(t_k^{2} - c_k\kappa - c_k)\cdot \frac{(r_{k} + 2)^{2}}{2\sigma^{2}}\right)},
\end{align*}
where $X\sim \mathbb{P}_{n_k}*\mathcal{N}(0, \sigma^2)$. Hence we obtain that
$$\mathbb{E}[W_{2}(\mathbb{P} * \mathcal{N}(0, \sigma^2), \mathbb{P}_{n_k} * \mathcal{N}(0, \sigma^2))]\ge \frac{\sqrt{C_l}}{64}\sqrt{\exp\left(-(t_k^{2} - c_k\kappa - c_k)\cdot \frac{(r_{k} + 2)^{2}}{2\sigma^{2}}\right)},$$
which indicates that there exists some constant $C_5, C_6$ only depending on $C, K, \sigma$ such that
$$\frac{\mathbb{E}[W_{2}(\mathbb{P} * \mathcal{N}(0, \sigma^2), \mathbb{P}_{n_k} * \mathcal{N}(0, \sigma^2))]}{n_k^{-\frac{(t^{2} - c\kappa - c)/(4\sigma^{2})}{(t^{2} - c\kappa - c)/\sigma^{2} - c^{2}/(2K^{2})}}}\ge C_5\exp\left(-C_6 r_k\right)\ge n_k^{-\mathcal{O}\left(\frac{1}{\sqrt{\log n_k}}\right)}.$$
Next we remember that $t_k = \frac{1}{2}(1 + \kappa)(1 + c_k)$, therefore if choosing $c$ large enough, we will have
$$t_k^2 - c_k\kappa - c_k = \frac{(1 + \kappa)^2(1+c_k)^2}{4} - c_k(1 + \kappa) = \frac{(1 + \kappa)^2c_k^2}{4} + \mathcal{O}(c_k),$$
which indicates that
$$\frac{(t^{2} - c\kappa - c)/(4\sigma^{2})}{(t^{2} - c\kappa - c)/\sigma^{2} - c^{2}/(2K^{2})} = \frac{(1 + \kappa)^2c_k^2 + \mathcal{O}(c_k)}{4(1 + \kappa)^2c_k^2 + \mathcal{O}(c_k) - 8c_k^2\kappa} = \frac{(1 + \kappa)^2}{4(1 + \kappa)^2 - 8\kappa} + \mathcal{O}\left(\frac{1}{c_k}\right) = \alpha + \mathcal{O}\left(\frac{1}{c_k}\right).$$
Therefore, choosing $c_k = M^k$ with $M = \max\left\{\sqrt{\frac{2}{\kappa}}, \frac{\kappa + 3}{1 - \kappa}, 3\right\}$, then for every $k$ we have $c_k\ge \max\left\{\sqrt{\frac{2}{\kappa}}, \frac{\kappa + 3}{1 - \kappa}, 3\right\}$, which indicates that this choice of $c_k$ satisfies all previous assumptions on $c_k$. We further notice that $r_k = M^{\frac{k(k-1)}{2}}$, hence
$$\frac{n_k^{-\frac{(t^{2} - c\kappa - c)/(4\sigma^{2})}{(t^{2} - c\kappa - c)/\sigma^{2} - c^{2}/(2K^{2})}}}{n_k^{-\alpha}}\ge n_k^{-\mathcal{O}\left(\frac{1}{\sqrt{\log \log n_k}}\right)}.$$
Therefore, we obtain that
$$n^\alpha\mathbb{E}[W_{2}(\mathbb{P} * \mathcal{N}(0, \sigma^2), \mathbb{P}_{n_k} * \mathcal{N}(0, \sigma^2))]\le n_k^{-\mathcal{O}\left(\frac{1}{\sqrt{\log \log n_k}}\right)}.$$
We let $k$ goes to infinity, and obtain that
$$\mathop{\limsup}_{n\to\infty}n^{\alpha + \frac{1}{\sqrt[3]{\log\log n}}}\mathbb{E}[W_{2}(\mathbb{P} * \mathcal{N}(0, \sigma^2), \mathbb{P}_{n} * \mathcal{N}(0, \sigma^2))]\ge \limsup_{k\to\infty} n_k^{\frac{1}{\sqrt[3]{\log\log n_k}} - \mathcal{O}\left(\frac{1}{\sqrt{\log \log n_k}}\right)} > 0.$$
And the proof of the lower bound part of Theorem \ref{thm-1D} is completed.
\end{proof}

\section{Proof of the Upper Bound of Theorem \ref{thm-1D}}\label{thm-upper-bound-part}

Without loss of generality, we consider the case $\sigma = 1$, as we can always reduce to this by rescaling.  We start the proof from the following observation~\citep[Theorem 2.18]{villani2003topics}: 
for two distributions $\mathbb{Q}_1, \mathbb{Q}_2$ on $\mathbb{R}$ with absolutely continuous CDFs
$F_{1}(x), F_{2}(x)$ the optimal coupling for the 2-Wasserstein distance is given by
$F_2^{-1}(F_1(x))$, implying an explicit formula:
	$$W_{2}(\mathbb{Q}_1, \mathbb{Q}_2)^{2} =\EE[\left|F_{2}^{-1}(F_{1}(T_1)) -
	T_1\right|^{2}]\,, T_1
	\sim \mathbb{Q}_1\,,$$
	where $F_{2}^{-1}(\cdot)$ is the inverse function of $F_{2}(\cdot)$. 
In our case we set $\mathbb{Q}_1 = \mathbb{P}*\calN(0,1)$ and $\mathbb{Q}_2 = \PP_n * \calN(0,1)$, and denote their CDFs
by $F$ and $F_n$ respectively. In the following, we use $\mathcal{N}$ to denote $\mathcal{N}(0, 1)$. We also denote by $\rho(t)$ the pdf of $F$ and by $T$ the optimal
transport map
	$$ T(t) \eqdef F_n^{-1}(F(t))\,.$$
Note that because of the randomness of $F_n$ the map $T$ is random as well. Our proof will proceed
along the following reductions:
\begin{enumerate}
\item (Conditioning) Note that if $E$ is any event with probability at least $1-O({1\over n^2})$ then we have
	\begin{equation}\label{eq-conditioning}\EE[W_2^2(\PP*\calN,\PP_n*\calN)] \le \EE[W_2^2(\PP*\calN,\PP_n*\calN)|E] +
	O({1\over n})\end{equation}
	This allows us to condition on a typical realization of the empirical measure $\PP_n$.


\item (Truncation) We will show that with high probability
	\begin{equation}\label{eq:pcp_2}
		|T(t) - t| \lesssim |t| + \sqrt{\log n} \qquad \forall t\in \mreals\,.
\end{equation}	
	Conditioning on this event, then, allows us to restrict evaluation of $W_2^2$ to a $O(\log n)$ range of $t$:
	$$ W_2^2(\PP*\calN,\PP_n*\calN) = \EE[|T(X+Z) - (X+Z)|^2 1\{|X+Z| \le b\sqrt{\log n}] +
	O(1/n)\,.$$

\item (Key bound) So far we are left to bound the integral
	\begin{equation}\label{eq:pcp_4}
		\int_{|t| \le b\sqrt{\log n}} \rho(t) |T(t)-t|^2 dt
\end{equation}	
	and we only have the bound~\eqref{eq:pcp_2}. The \textit{key novel ingredient} is the
	following observation. The transport distance can be
	bounded by 
	\begin{equation}\label{eq:pcp_3}
		|T(t) - t| \le {|F(t) - F_n(t)|\over \rho(t)}.
	\end{equation}	
	This bound can be explained by the fact that if $F(t) < F_n(t)$ then the distance we need to
	travel to the right of $t$ so that $F(\cdot)$ raise to the value of $F_n(t)$ will be
	around $|F(t) - F_n(t)|\over \rho(t)$. There are several caveats in a rigorous statement
	of~\eqref{eq:pcp_3} (see Prop.~\ref{prop-cdf-pdf} for details), but the most important one
	is that it only holds provided the RHS of~\eqref{eq:pcp_3} is $\le 1$. 
\item (Concentration) Next, we will show that with high probability
	\begin{equation}\label{eq:pcp_5}
		|F_n(t) - F(t)| \lesssim {\log n\over \sqrt{n}} \sqrt{\min(F(t),1-F(t)) \vee {1\over n}}\,. 
\end{equation}
	It turns out that we also have $\min(F(t),1-F(t)) \lesssim \rho(t)^{4K^2\over (1+K^2)^2} n^{o(1)}$.
	Hence, we have a transport distance bound
		$$ |T(t) - t| \lesssim {\log n \over \sqrt{n}} \rho(t)^{{2K^2\over (1+K^2)^2} - 1}
		n^{o(1)}\,,
		$$
		provided the RHS is $\lesssim 1$, which is equivalent to say $\rho(t) >
		n^{-\alpha - o(1)}$ for $\alpha = \frac{(1 + K^2)^2}{2(1 + K^4)}$. Note that in this region 
		the integral~\eqref{eq:pcp_4} becomes bounded by 
			\begin{equation}\label{eq:pcp_7}
				n^{o(1)} \int_{|t| \le b\sqrt{\log n}, \rho(t)>n^{-\alpha}} {1\over n}
			\rho^{\frac{4K^2}{(1 + K^2)^2} - 1}(t) \le n^{-\alpha + o(1)}\,,
\end{equation}			
			since $K>1$ and thus the power of $\rho$ is negative.
\item (Final) The final step is to split the integral~\eqref{eq:pcp_4} into values of $\rho(t)<
n^{-\alpha}$ (for which we use the bound~\eqref{eq:pcp_2} and $|t| \lesssim \sqrt{\log n}$) and $\rho(t) > n^{-\alpha}$ (for which
we use~\eqref{eq:pcp_7}). This gives us the contributions
	$$ n^{-\alpha} O(\log n) + n^{o(1)} n^{-\alpha}$$
	completing the proof.
\end{enumerate}

\par If we can prove for any $\epsilon > 0$, there exists $C_\epsilon$ such that for any $n$ and $K$-subgaussian distribution $\PP$,
\begin{equation}\label{eq-upper-bound1}\EE[W_2^2(\PP * \mathcal{N}, \PP_n * \mathcal{N})]\le C_{\epsilon} n^{-2\alpha + \epsilon},\end{equation}
then for every integer $t$ and $n$ we have
$$\EE[W_2^2(\PP * \mathcal{N}, \PP_n * \mathcal{N})]\le C_{1/(2t)} n^{-2\alpha + 1/(2t)}.$$
WLOG we assume that $C_{1/(2t)}\ge C_{1/(2s)}\ge 1$ for every $t > s$. Therefore, when $n\ge C_{1/(2t)}^{2t}$ we have $\EE[W_2^2(\PP * \mathcal{N}, \PP_n * \mathcal{N})]\le n^{-2\alpha + 1/t}$ for all $K$-subgaussian distribution $\PP$. We let $\delta_n = 1/(2t)$ for those $n\in (C_{1/(2(t-1))}^{2(t-1)}, C_{1/(2t)}^{2t}]$, and for those $n\le C_{1/2}^2$, we choose $\delta_n = \log_2\left[\max_{2\le n\le C_{1/2}^2}\EE[W_2^2(\PP * \mathcal{N}, \PP_n * \mathcal{N})]\right]$, we will have
$$\EE[W_2^2(\PP * \mathcal{N}, \PP_n * \mathcal{N})]\le n^{-2\alpha + \delta_n}$$
with $\lim_{n\to\infty}\delta_n = 0$. Therefore, we only need to prove \eqref{eq-upper-bound1}

\begin{proposition}\label{prop-PDF-CDF}
    We denote the CDF, PDF of $\mathbb{P} * \mathcal{N}(0, 1)$ as $F, \rho$, respectively, and let
    $X\sim \PP$. Suppose there exist constants $C, K > 0$ such that for $\forall r\ge 0$,
    \begin{equation}\label{eq:proposition-subgaussian}\mathbf{P}(|X|\ge r)\le C\exp\left(-\frac{r^2}{2K^2}\right).\end{equation}
   	For $\beta = \frac{4K^2}{(1 + K^2)^2}$ and any $0 < \epsilon < \beta$, $\exists M = M(K, C, \epsilon)\ge 1$ such that for any $K$-subgaussian distribution $\PP$,
    \begin{align*}
    	1 - F(r) \le M\rho(r)^{\beta - \epsilon}, \quad \forall r\ge 0,\qquad\text{and}\qquad F(r) \le M\rho(r)^{\beta - \epsilon}, \quad \forall r < 0.
	\end{align*}
\end{proposition}
\begin{remark}
	We notice that this result is tight when considering $\PP_h = (1-p_h)\delta_0 + p_h\delta_h$ with $p_h = \exp\left(-\frac{h^2}{2K^2}\right)$ and $r = \frac{(K^2 + 1)h}{2K^2}$. Then we have $\rho(r) = \Theta\left(\exp\left(-\frac{(K^2+1)^2h^2}{8K^4}\right)\right)$ and $1 - F(r) = \Theta\left(\exp\left(-\frac{h^2}{2K^2}\right)\right)$ when $h\to \infty$. Hence the above inequalities are tight.
\end{remark}


First we present two lemmas:
\begin{lemma}\label{lem-cdf-pdf}
	Suppose $\Phi_{1}$ to be the CDF of Gaussian distribution $\mathcal{N}(0, 1)$, then we have $1 - \Phi_{1}(l) \le \exp(-l^2/2)$ for any $l\le 0$ and $\Phi_{1}(l) \le \exp(-l^2/2)$ for any $l < 0$.
\end{lemma}
\begin{proof}
	This lemma directly follows from \citep[Proposition 2.1.2]{vershynin2018high}.
\end{proof}
\begin{lemma}\label{lem-cr}
	$\mathbb{P}$ is a 1-dimensional $K$-subgaussian distribution (For $X\sim \mathbb{P}$ and any $r\ge 0$, $\mathbf{P}(X\ge r)\le C\exp\left(-\frac{r^{2}}{2K^{2}}\right)$ for constant $C > 0$), and $\rho(\cdot)$ is the PDF of $\mathbb{P} * \mathcal{N}(0, 1)$. For any $0 < \epsilon < \beta$ with $\beta = \frac{4K^2}{(1 + K^2)^2}$ we have
	$$\rho(r)\ge C_{\epsilon}\mathbf{P}(X\ge r)^{\frac{1}{\beta - \epsilon}}\quad \forall r\ge 0$$
	for some positive constant $C_{\epsilon} = C_\epsilon(K, C)$.
\end{lemma}
\begin{proof}
		\par For $X\sim\mathbb{P}$, choosing $M = M(K, C)\triangleq K\sqrt{2\log (2C)}> 0$ then we have
		$$\mathbf{P}(X\in [-M, M]) = 1 - \mathbf{P}(|X| > M)\le 1 - C\exp\left(-\frac{M^2}{2K^2}\right)\ge\frac{1}{2}.$$
		For $0\le r\le M$, we have 
	$$\rho(r)\ge \int_{-M}^{M}\eta(x)\varphi_{1}(r-x)dx\ge \mathbf{P}(X\in [-M, M])\cdot \min_{-M\le x\le M}\varphi(r - x)\ge \frac{1}{2}\varphi_{1}(2M),$$
	where we use $\eta(\cdot)$ to denote the PDF of $\mathbb{P}$ (which can be a generalized function). Hence $\rho(r)\ge C_{\epsilon}\mathbf{P}(X\ge r)^{\frac{1}{\beta - \epsilon}}$ holds for all $0\le r\le M$ if $C_{\epsilon}\le \frac{1}{2}\varphi_1(2M)$.
	\par Next we consider cases where $r\ge M$. We let $c_{r} = \log\frac{C}{\mathbf{P}(X\ge r)}$. If $c_r\ge \log\frac{1}{\rho(r)}$, then we have 
	$$\rho(r)\ge \frac{\mathbf{P}(X\ge r)}{C}\ge \frac{1}{C}\mathbf{P}(X\ge r)^{\frac{1}{\beta - \epsilon}},$$
	where we use the fact that $\beta - \epsilon\le \beta\le 1$ and hence $\frac{1}{\beta - \epsilon}\ge 1$. 
	
	\par Next we consider cases where $c_r < \log\frac{1}{\rho(r)}$. We let $r_{1} = \sqrt{2\log (2) + 2c_{r}}K$, then we have $\mathbf{P}(X\ge r) = Ce^{-c_{r}}$, and
	$$\mathbf{P}(X\ge r_{1})\le C\exp\left(-\frac{r_{1}^{2}}{2K^{2}}\right)\le \frac{C}{2}e^{-c_{r}}.$$
	Hence $\mathbf{P}(r < X\le r_{1})\ge \frac{C}{2}e^{-c_{r}}$, which indicates that
	\begin{align*}
		\rho(r) & \ge \int_{-M}^{M}\eta(x)\varphi(r-x)dx + \int_{r}^{r_{1}}\eta(x)\varphi(r-x)dx\ge \frac{1}{2}\varphi(r+M) + \frac{C}{2}e^{-c_{r}}\varphi(r_{1}-r)\\
	& = \frac{1}{2\sqrt{2\pi}}\exp\left(-\frac{(r+M)^{2}}{2}\right) + \frac{C}{2\sqrt{2\pi}}\exp\left(-c_{r} - \frac{(r_{1}-r)^{2}}{2}\right).
	\end{align*}
	Next, we let $c_r = \frac{x^2r^2}{2K^2}$ with $x\ge 1$. We notice that when $x\ge \frac{2K^2}{1 + K^2}$ we have $\frac{x^2}{\beta K^2}\ge 1$, hence $-\frac{r^2}{2}\ge -\frac{1}{\beta}c_r$; and when $1\le x\le \frac{2K^2}{1 + K^2}$, we have $\left(\frac{1}{\beta K^2} - \frac{1}{K^2} - 1\right)x^2 + 2x - 1\ge 0$ since $-1\le \frac{1}{\beta K^2} - \frac{1}{K^2} - 1 = \frac{(K^2+1)(1-3K^2)}{4K^4}\le 0$, and hence $-c_r - \frac{(r - \sqrt{2c_r}K)^2}{2}\ge -\frac{1}{\beta}c_r$. Therefore, we have
	$$\max\left\{-\frac{r^2}{2}, -c_r - \frac{(r - \sqrt{2c_r}K)^2}{2}\right\}\ge -\frac{1}{\beta} c_r$$
	
	We further notice that
	\begin{align*}
		\rho(r) & = \int_{-\infty}^\infty\eta(x)\varphi_1(r - x)dx = \int_{-\infty}^{r/2}\eta(x)\varphi_1(r-x)dx + \int_{r/2}^{\infty}\eta(x)\varphi_1(r-x)dx\\
		& \le \sup_{x\le r/2}\varphi_1(r-x) + \mathbf{P}(X\ge r/2]\le \frac{1}{\sqrt{2\pi}}\exp\left(-\frac{r^2}{8}\right) + C\exp\left(-\frac{r^2}{8K^2}\right)\le \left(\frac{1}{\sqrt{2\pi}} + C\right)\exp\left(-\frac{r^2}{8K^2}\right). \numberthis \label{eq-form-rho}
	\end{align*}
	This indicates that $r\le 2\sqrt{K\log \frac{\bar{C}}{\rho(r)}}$ with $\bar{C} = \frac{1}{\sqrt{2\pi}} + C$, and hence $\exp(rM) = \mathcal{O}(\rho(r)^{-\epsilon'})$ for $\forall\epsilon' > 0$. We further notice that $\exp\left(2\sqrt{\log(2)}K^2\sqrt{c_r}\right) = \mathcal{O}(\rho(r)^{-\epsilon'})$ for $\forall\epsilon' > 0$ since $c_r\le \log\frac{1}{\rho(r)}$. Therefore, we obtain that
	\begin{align*}
		&\quad \mathbf{P}(X\ge r)^{\frac{1}{\beta}} = C^{(1/\beta)} \exp\left(-\frac{1}{\beta}c_r\right)\\
		& \le C^{(1/\beta)}\max\left\{\exp\left(-\frac{r^2}{2}\right), \exp\left(-c_r - \frac{(r - \sqrt{2c_r}K)^2}{2}\right)\right\}\\
		&\le \frac{2\sqrt{2\pi}C^{(1/\beta)}}{1 + C}\cdot \left(\frac{1}{2\sqrt{2\pi}}\exp\left(-\frac{r^2}{2}\right) + \frac{C}{2\sqrt{2\pi}}\exp\left(-c_r - \frac{(r - \sqrt{2c_r}K)^2}{2}\right)\right)\\
		&\le \rho(r)\cdot \max\left\{\exp\left(Mr + \frac{M^2}{2}\right), \exp\left(\sqrt{2\log(2)}K(\sqrt{2c_r}K - r) + 2K^2\log 2\right)\right\}\\
		&\le \rho(r)\cdot \max\left\{\exp\left(Mr + \frac{M^2}{2}\right), \exp\left(2\sqrt{\log(2)}K^2\sqrt{c_r} + 2K^2\log 2\right)\right\}\\
		&\le \rho(r)\cdot \tilde{\mathcal{O}}(\rho^{-\epsilon'}).
	\end{align*}
	Choosing $\epsilon' = \frac{\epsilon}{\beta}$, we know that there exists some positive constant $C_{\epsilon}$ such that
	$$\rho(r)\ge C_{\epsilon}\mathbf{P}(X\ge r)^{\frac{1}{\beta - \epsilon}}\quad \forall r\ge 0.$$
\end{proof}

\begin{proof}[Proof of Proposition \ref{prop-PDF-CDF}]
    We only prove this results for $r\ge 0$, as the proof of $r\le 0$ is similar. First we can write
    \begin{align*}
		1 - F(r) = \int_{-\infty}^{\infty}\eta(t)(1 - \Phi_1(r - t))dt,\quad \text{and}\quad \rho(r) = \int_{-\infty}^{\infty}\eta(t)\cdot\frac{1}{\sqrt{2\pi}}\exp\left(-\frac{(r - t)^2}{2}\right)dt.
    \end{align*}
	Noticing that $\mathbf{P}(|X|\ge r)\le C\exp\left(-\frac{r^2}{2K^2}\right)$, with $\tilde{K}\triangleq K\sqrt{2(\log(2C))}$ we have $\mathbf{P}(|X|\ge \tilde{K})\le C\exp(-\log (2C)) = 1/2$ and hence $\mathbf{P}(|X|\le \tilde{K})\ge 1/2$. In the following, we will discuss cases where $0\le r\le \tilde{K}$ and $r > \tilde{K}$ separately.
    \par If $0\le r\le \tilde{K}$, then we have
    \begin{align*}
    	\sqrt{2\pi}\rho(r) \ge \int_{-\tilde{K}}^{\tilde{K}}\rho(t)\exp\left(-\frac{(r - t)^2}{2}\right)dt\ge \mathbf{P}(|X|\le \tilde{K})\cdot \min_{\substack{0\le r\le \tilde{K}\\t\in [-\tilde{K}, \tilde{K}]}}\exp\left(-\frac{(r - t)^2}{2}\right) = \frac{1}{2}\exp\left(-2\tilde{K}^{2}\right),
	\end{align*}
	which implies that with any $\varepsilon > 0$ and some constant $M_1 = M_1(\varepsilon) > 0$, $1 - F(r)\le 1\le M_{1}\rho(r)^{\beta - \epsilon}$ for any $r\in [0, R_{0}]$.
    
    \par Next consider cases where $r > \tilde{K}$. According to \eqref{eq:proposition-subgaussian}, we have $\mathbf{P}(X \ge r)\le C\exp(-r^2/(2K^2))$, which implies that
    \begin{align*}
    	1 - F(r) & = \int_{-\infty}^{r}\eta(t)(1 - \Phi_1(r - t))dt + \int_{r}^{\infty}\eta(t)(1 - \Phi_1(r - t))dt\le \int_{-\infty}^{r}\eta(t)(1 - \Phi_1(r - t))dt + \mathbf{P}(X \ge r),
	\end{align*}
	where in the last inequality we use the fact that $\mathbf{P}(X \ge r) = \int_{-\infty}^{r}\eta(t)dt$. For $r > t$, according to Lemma \ref{lem-cdf-pdf}, we have $1 - \Phi_1(r - t)\le \exp(-(r - t)^2/2)$. Hence we have
	\begin{align*}
		1 - F(r) & \le \int_{-\infty}^{r}\eta(t)\exp\left(-\frac{(r - t)^2}{2}\right)dt + \mathbf{P}(X \ge r) \le \sqrt{2\pi}\cdot \rho(r) + \mathbf{P}(X \ge r),
	\end{align*}
	where in the last inequality we use $\rho(r) = \int_{-\infty}^{\infty}\eta(t)\exp(-(r - t)^2/2)dt$. Moreover, according to Lemma \ref{lem-cr}, there exists a constant $C_{\epsilon}$ such that $\rho(r)\ge C_{\epsilon}\mathbf{P}(X\ge r)^{\frac{1}{\beta - \epsilon}}$, which indicates that $\mathbf{P}(X \ge r)\le C_{\epsilon}^{-\beta + \epsilon}\rho(r)^{\beta - \epsilon}$. Hence we have
	$$1 - F(r)\le \sqrt{2\pi}\cdot\rho(r) + C_{\epsilon}^{-\beta + \epsilon}\rho(r)^{\beta - \epsilon}.$$
    When $\rho(r)\le 1$, $\rho(r)\le \rho(r)^{\beta - \epsilon}$. Hence $1 - F(r)\le (C_{\epsilon}^{-\beta + \epsilon} + \sqrt{2\pi})\cdot \rho(r)^{\beta - \epsilon}.$ When $\rho(r) > 1$, $\rho(r)^{\beta - \epsilon} > 1$. Hence $1 - F(r)\le 1 < \rho(r)^{\beta - \epsilon}$. Therefore, by choosing $M = \max\{M_{1}, (C_{\epsilon}^{-\beta + \epsilon} + \sqrt{2\pi}), 1\}\ge 1$, we obtain
    $$1 - F(r)\le M\rho(r)^{\beta - \epsilon}, \quad \forall r\ge 0.$$
\end{proof}

\begin{proposition}\label{prop-cdf-pdf}
    Consider two distributions $\mathbb{P}, \mathbb{Q}$ on $\mreals$. We denote the PDF and CDF of
    $\mathbb{P}$ as $\rho_\PP$ and $F_\PP$; and similarly for $\QQ$. Assume that
    $\rho_\QQ>0$ everywhere on $\mreals$. 
    For a fixed $h>0$ denote $L_h(t) \eqdef \sup_{x\in [t - h, t +
    h]}|F_\PP(x) - F_\QQ(x)|$ and $\underline{\rho}_h(t) = \inf_{x\in [t - h, t + h]}\rho_\PP(x)$. If
    we have 
	$$\Delta_{h}(t) \eqdef {L_h(t) \over \underline{\rho}_h(t)} \le h,$$
	then 
	$$\left|F_\QQ^{-1}(F_\PP(t)) - t\right|\le \Delta_h(t)\,. $$
    \end{proposition}

\begin{proof}
	Suppose that $F_{\PP}(t) \ge F_{\QQ}(t)$ and let $h' = \Delta_h(t) \le h$. Then we claim that
	\begin{equation}\label{eq:pcp_1}
		F_{\PP}(t) \le F_{\QQ}(t+h')\,.  
\end{equation}	
	Indeed, we have $F_{\PP}(t) \le F_{\PP}(t+h') - \underline{\rho}_h(t) h'  =
	F_{\PP}(t+h') - L_h(t)\le F_{\PP}(t+h') - (F_{\PP}(t+h') - F_{\QQ}(t+h')) = F_{\QQ}(t+h')$. Now, since $F_{\QQ}(t) \le F_{\PP}(t)\le F_{\QQ}(t+h')$ we obtain that
	$0 < F_{\QQ}^{-1}(F_{\PP}(t)) - t \le h'$. The case of $F_{\PP}(t)<F_{\QQ}(t)$ is treated similarly.
\end{proof}

\begin{proposition}\label{prop-pq}
	Assume distribution $\PP, \QQ$ are $K_1, K_2$-subgaussian distributions respectively, e.g. for any $X\sim\PP, Y\sim \QQ$ we have 
	$$\mathbf{P}(|X|\ge r)\le C_{1}\exp\left(-\frac{r^{2}}{2K_{1}^{2}}\right), \quad \mathbf{P}(|Y|\ge r)\le C_{2}\exp\left(-\frac{r^{2}}{2K_{2}^{2}}\right), \quad \forall r>0\,.$$
	Then for all $x\in\mathbb{R}$, we have
	$$\left|F_{\QQ, 1}^{-1}(F_{\PP, 1}(x)) - x\right|\le 2|x| + 2 + \tilde{K}_{1} + \tilde{K}_{2}(|x| + 2 + \tilde{K}_{1})\,,$$
	where $F_{\PP, 1}, F_{\QQ, 1}$ denote the CDFs of $X+Z$ and $Y+Z$, $(X,Y)\dperp Z \sim\mathcal{N}$, and $\tilde{K}_{1}\triangleq K_{1}\sqrt{2\log 2C_{1}}, \tilde{K}_{2}(t)\triangleq K_{2}t + K_{2}\sqrt{2\log\left(4tC_{2}\right)}$.
	\end{proposition}
\begin{proof}
	We let $R = |x|$ and $\tilde{R} = R + 2 + \tilde{K}_{1}$. First we notice that the PDFs of distribution $\mathbb{P} * \mathcal{N}, \mathbb{Q} * \mathcal{N}$ at any real number is positive, hence $F_{\PP, 1}, F_{\QQ, 1}$ are monotonically increasing in the entire real line. We have $\mathbf{P}\left(|X|\ge \tilde{K}_{1}\right)\le C_{1}\exp\left(-\log(2C_{1})\right) = \frac{1}{2}.$ Therefore, we obtain that $\mathbf{P}\left(|X|\le \tilde{K}_{1}\right)\ge 1 - 1/2 = 1/2.$ We further notice that if $X\sim \mathbb{P}, Z\sim \mathcal{N}(0, 1)$ are independent, $X + Z\sim \mathbb{P} * \mathcal{N}(0, 1)$. And also
	\begin{align*}
		\{|X|\le \tilde{K}_{1}\}\cap\{Z\le -\tilde{K}_{1} - R\} & \subset \left\{X + Z \le -R\right\}\\
		\{|X|\le \tilde{K}_{1}\}\cap\{Z\ge \tilde{K}_{1} + R\} & \subset \left\{X + Z \ge R\right\}.
	\end{align*}
	Recall that we use $\Phi_{1}$ to denote the CDF of distribution $\mathcal{N}(0, 1)$. Hence noticing that $\Phi_{1}(-R - \tilde{K}_{1}) = 1 - \Phi_{1}(R + \tilde{K}_{1}) = \mathbf{P}(Z \le -\tilde{K}_{1} - R) = \mathbb{P}(Z\ge \tilde{K}_{1} + R)$, we have 
	\begin{align*}
		\frac{1}{2}\Phi_{1}(-R-\tilde{K}_{1}) & \le \mathbf{P}\left(|X|\ge \tilde{K}_{1}\right)\mathbf{P}(Z\le -\tilde{K}_{1} - R)\le \mathbf{P}(X + Z\le -R) = F_{\PP, 1}(-R)\\
		\frac{1}{2}\Phi_{1}(-R-\tilde{K}_{1}) & \le \mathbf{P}\left(|X|\ge \tilde{K}_{1}\right)\mathbf{P}(Z\ge \tilde{K}_{1} + R)\le \mathbf{P}(X + Z \ge R) = 1 - F_{\PP, 1}(R),
	\end{align*}
	which indicates that
	$$\frac{1}{2}\Phi_{1}(-R-\tilde{K}_{1})\le F_{\PP, 1}(-R)\le F_{\PP, 1}(R)\le 1 - \frac{1}{2}\Phi_{1}(-R-\tilde{K}_{1}).$$
	\par Next, if $Y\sim \mathbb{Q}, Z\sim \mathcal{N}(0, 1)$ are independent, we have $Y + Z\sim \mathbb{Q} * \mathcal{N}(0, 1)$. Noticing that,
	\begin{align*}
		\left\{Y + Z \le - \tilde{R} - \tilde{K}_{2}(\tilde{R})\right\} & \subset\{Z\le - \tilde{R}\}\cup \{Y\le -\tilde{K}_{2}(\tilde{R})\},\\
		\left\{Y + Z \ge \tilde{R} + \tilde{K}_{2}(\tilde{R})\right\} & \subset\{Z\ge \tilde{R}\}\cup \{Y\ge \tilde{K}_{2}(\tilde{R})\},
	\end{align*}
	we obtain that
	\begin{align*}
		F_{\QQ, 1}(-\tilde{R} - \tilde{K}_{2}(\tilde{R})) & \le \Phi_{1}(- \tilde{R}) + \mathbf{P}(|Y|\ge \tilde{K}_{2}(\tilde{R})),\\
		1 - F_{\QQ, 1}(\tilde{R} + \tilde{K}_{2}(\tilde{R})) & \le 1 - \Phi_{1}(\tilde{R}) + \mathbf{P}(|Y|\ge \tilde{K}_{2}(\tilde{R})) = \Phi_{1}(-\tilde{R}) + \mathbf{P}(|Y|\ge \tilde{K}_{2}(\tilde{R})).
	\end{align*}
	According to Proposition 2.1.2 in \cite{vershynin2018high}, we have
	$$\Phi_{1}(- \tilde{R})\ge\left(\frac{1}{\tilde{R}} - \frac{1}{\tilde{R}^{3}}\right)\cdot \frac{1}{\sqrt{2\pi}}\exp\left(-\frac{\tilde{R}^{2}}{2}\right).$$
	Hence since $\tilde{R} = R + \tilde{K}_{1} + 2\ge 2$, we will have
	$$\Phi_{1}( - \tilde{R})\ge \frac{3}{4\tilde{R}}\cdot \frac{1}{\sqrt{2\pi}}\exp\left(-\frac{\tilde{R}^{2}}{2}\right)\ge \frac{1}{4\tilde{R}}\exp\left(-\frac{\tilde{R}^{2}}{2}\right).$$
	We further notice that $\tilde{K}_{2}(\tilde{R}) = K_{2}\tilde{R} + K_{2}\sqrt{2\log\left(4\tilde{R}C_{2}\right)}$, hence we have 
	$$\mathbf{P}(|Y|\ge \tilde{K}_{2}(\tilde{R}))\le C_{2}\exp\left(-\frac{\tilde{K}_{2}(\tilde{R})^{2}}{2K_{2}^{2}}\right) = \frac{1}{4\tilde{R}}\exp\left(-\frac{\tilde{R}^{2}}{2}\right)\le \Phi_{1}( - \tilde{R}),$$
	which indicates that 
	$$F_{\QQ, 1}(-\tilde{R} - \tilde{K}_{2}(\tilde{R}))\le 2\Phi_{1}(- \tilde{R}) \quad\text{and}\quad 1 - F_{\QQ, 1}(\tilde{R} + \tilde{K}_{2}(\tilde{R}))\le 2\Phi_{1}(- \tilde{R}).$$
	
	\par Additionally, since for any $t\le 0$, we have
	$$\exp\left(-\frac{(t - 2)^{2}}{2}\right)\le \exp\left(-\frac{t^{2}}{2} - \frac{4}{2}\right) = \exp(-2)\cdot \exp\left(-\frac{t^{2}}{2}\right)\le \frac{1}{4}\exp\left(-\frac{t^{2}}{2}\right).$$
	This indicates that
	\begin{align*}
		\frac{1}{4}\Phi_{1}(-R- \tilde{K}_{1}) & = \frac{1}{4}\cdot \frac{1}{\sqrt{2\pi}}\int_{-\infty}^{-R - \tilde{K}_{1}}\exp\left(-\frac{t^{2}}{2}\right)dt \ge \frac{1}{\sqrt{2\pi}}\int_{-\infty}^{-R - \tilde{K}_{1}}\exp\left(-\frac{(t - 2)^{2}}{2}\right)dt\\
		& = \frac{1}{\sqrt{2\pi}}\int_{-\infty}^{-R - \tilde{K}_{1} - 2}\exp\left(-\frac{t^{2}}{2}\right)dt = \Phi_{1}(-R- \tilde{K}_{1} - 2).
	\end{align*}
	Therefore, we obtain that
	$$F_{\QQ, 1}\left(-\tilde{R} - \tilde{K}_{2}(\tilde{R})\right)\le 2\Phi_{1}\left(-\tilde{R}\right) = 2\Phi_{1}(-R - \tilde{K}_{1}-2)\le \frac{1}{2}\Phi_{1}(-R - \tilde{K}_{1})\le F_{\PP, 1}(-R).$$
	Similarly, we can also obtain that $F_{\PP, 1}(R)\le F_{\QQ, 1}\left(\tilde{R} + \tilde{K}_{2}(\tilde{R})\right)$. Hence using the monotonicity of $F_{\PP, 1}$ and $F_{\QQ, 1}$, we obtain that, 
	$$F_{\QQ, 1}\left(-\tilde{R} - \tilde{K}_{2}(\tilde{R})\right)\le F_{\PP, 1}(-R)\le F_{\PP, 1}(x)\le F_{\PP, 1}(R)\le F_{\QQ, 1}\left(\tilde{R} + \tilde{K}_{2}(\tilde{R})\right),$$
	which indicates that
	$$ - \tilde{R} - \tilde{K}_{2}(\tilde{R})\le F_{\QQ, 1}^{-1}(F_{\PP, 1}(x))\le \tilde{R} + \tilde{K}_{2}(\tilde{R}).$$
	Hence we have
	$$\left|F_{\QQ, 1}^{-1}(F_{\PP, 1}(x)) - x\right|\le R + \tilde{R} + \tilde{K}_{2}(\tilde{R}) = 2|x| + 2 + \tilde{K}_{1} + \tilde{K}_{2}(|x| + \tilde{K}_{1} + 2).$$
\end{proof}

\begin{proposition}\label{prop-difference}
    Suppose $F, F_{n}$ are 
    CDFs of distribution $\mathbb{P}*\mathcal{N}$ and $\mathbb{P}_n*\mathcal{N}$ respectively. Then with probability at least $1 - \delta$, we have the following inequality:
    $$\sup_{t\in\mathbb{R}}\frac{|F(t) - F_{n}(t)|}{\sqrt{1/n\vee(F(t)\wedge (1-F(t)))}}\le \frac{16}{\sqrt{n}}\log\left(\frac{2n}{\delta}\right).$$
\end{proposition}
To prove this proposition, we first present a lemma indicating a similar result without Gaussian smoothing:
\begin{lemma}\label{lem-CDF}
    For a given distribution $\QQ$ on real numbers with always-positive PDF, we denote its empirical measure with $n$ data points to be $\QQ_n$ ($\QQ_n = \frac{1}{n}\sum_{i=1}^n\delta_{X_i}$ where $X_i\sim\QQ$ are $i.i.d.$). We further use $F, \hat{F}_n$ to denote the CDF of $\QQ, \QQ_n$ respectively.
    Then with probability at least $1 - \delta$, we have
    $$\sup_{t\in\mathbb{R}} \frac{|F(t) - \hat{F}_n(t)|}{\sqrt{1/n\vee(F(t)\wedge (1-F(t)))}}\le \frac{8}{\sqrt{n}}\log\left(\frac{n}{\delta}\right).$$
\end{lemma}
\begin{remark}
	From Theorem 2.1 of \cite{gine2006concentration} we can obtain a result similar to this lemma: if $\QQ$ is the uniform distribution on $[0, 1]$, then there exist universal positive constants $C_0, K$ such that for any $s > 0$,
	$$\mathbf{P}\left[\sup_{1/n\le t\le 1/2}\frac{|F(t) - \hat{F}_n(t)|}{\sqrt{t}}\ge \frac{4}{\sqrt{n}} + \frac{2s\log\log n}{\sqrt{n}\log\log\log n}\right]\le \frac{C_0}{\log(n)^{(s/(2K)-1)}}.$$
\end{remark}

\begin{remark}
    If we would like to obtain a uniform bound without truncation, then we have to pay an additional factor $\sqrt{1/\delta}$. This is summarized in the following results: with probability at least $1 - \delta$, we have
    $$\sup_{t\in\mathbb{R}} \frac{|F(t) - \hat{F}_n(t)|}{\sqrt{F(t)\wedge (1-F(t))}}\le 16\sqrt{\frac{1}{\delta n}}\log\left(\frac{4n}{\delta}\right).$$
    Also we have a lower bound to the LHS in the above inequality, indicating that the factor $\sqrt{1/\delta}$ is necessary: with probability at least $\delta$, we have
    $$\sup_{t\in\mathbb{R}} \frac{|F(t) - \hat{F}_n(t)|}{\sqrt{F(t)\wedge (1-F(t))}}\ge\sqrt{\frac{1}{2\delta n}}.$$
\end{remark}

\begin{proof}[Proof of Lemma \ref{lem-CDF}]
   	\par With loss of generality, we assume $\QQ$ is the uniform distribution on $[0, 1]$ (otherwise we consider the similar argument on distribution $\QQ(F^{-1}(\cdot))$). Then we have $F(t) = t$ for any $0\le t\le 1$. We only need to prove that with probability at least $1 - \delta$, 
    $$\sup_{t\in\mathbb{R}} \frac{|F(t) - \hat{F}_n(t)|}{\sqrt{1/n\vee(t\wedge (1-t))}}\le \sqrt{\frac{\log n}{n}}.$$
    
    According to Bernstein inequality, we have 
    $$\mathbf{P}\left(\left|F(t) - \hat{F}_n(t)\right|>\varepsilon\right)\le \exp\left(-\frac{n\varepsilon^2}{2t(1-t) + 2/3\varepsilon}\right)\le \exp\left(-\frac{n\varepsilon^2}{2t + 2/3\varepsilon}\right).$$
    Choosing $\varepsilon = 4\sqrt{\frac{t}{n}}\log\left(\frac{1}{\delta}\right)$, and noticing that with this choice we have $\frac{1}{2}n\varepsilon^2\ge 2t\log(1/\delta)$ and also $\frac{1}{2}n\varepsilon^2\ge \frac{2}{3}\varepsilon\log(1/\delta)$, we obtain that   
    $$\mathbf{P}\left(\left|F(t) - \hat{F}_n(t)\right|>
    4\sqrt{\frac{t}{n}}\log\left(\frac{1}{\delta}\right)\right)\le \delta.$$
    Therefore, choosing $t = \frac{k}{n}$ with $1\le k\le \frac{n}{2}$, and applying union bound for $1\le k\le \frac{n}{2}$, we obtain that
    $$\mathbf{P}\left(\left|F\left(\frac{k}{n}\right) - \hat{F}_n\left(\frac{k}{n}\right)\right|\le 4\sqrt{\frac{(k/n)}{n}}\log\left(\frac{n}{\delta}\right), \ \forall 1\le k\le \frac{n}{2}\right)\le \frac{\delta}{2}.$$
    We further notice that for any $\frac{k}{n}\le t\le \frac{k+1}{n}$, we have
    $$|F(t) - \hat{F}_n(t)| = |t - \hat{F}_n(t)|\le \frac{1}{n} + \max\left\{\left|F\left(\frac{k}{n}\right) - \hat{F}_n\left(\frac{k}{n}\right)\right|, \left|F\left(\frac{k+1}{n}\right) - \hat{F}_n\left(\frac{k+1}{n}\right)\right|\right\}.$$
    When $k\ge 1$ and $\frac{2k}{n}\le \frac{k+1}{n}$. Therefore, if for every $1\le k\le \frac{n}{2}$ we all have $\left|F\left(\frac{k}{n}\right) - \hat{F}_n\left(\frac{k}{n}\right)\right|\le 4\sqrt{\frac{(k/n)}{n}}\log\left(\frac{n}{\delta}\right)$, then for every $0\le t\le \frac{1}{n}$, we have
    $$\frac{|F(t) - \hat{F}_n(t)|}{\sqrt{1/n\vee(t\wedge (1-t))}}\le \frac{1/n + |F(1/n) - \hat{F}_n(1/n)|}{\sqrt{1/n}}\le 5\sqrt{\frac{1}{n}}\log\left(\frac{n}{\delta}\right),$$
    and for every $\frac{k}{n}\le t\le \frac{k+1}{n}$ with $k\le \frac{n}{2}$, we have
    \begin{align*}
    \frac{|F(t) - \hat{F}_n(t)|}{\sqrt{1/n\vee(t\wedge (1-t))}} & \le \frac{\frac{1}{n} + \max\left\{\left|F\left(\frac{k}{n}\right) - \hat{F}_n\left(\frac{k}{n}\right)\right|, \left|F\left(\frac{k+1}{n}\right) - \hat{F}_n\left(\frac{k+1}{n}\right)\right|\right\}}{\sqrt{k/n}}\\
    & \le \sqrt{\frac{1}{n}} + \sqrt{2}\cdot \max\left\{\frac{\left|F\left(\frac{k}{n}\right) - \hat{F}_n\left(\frac{k}{n}\right)\right|}{\sqrt{k/n}}, \frac{\left|F\left(\frac{k+1}{n}\right) - \hat{F}_n\left(\frac{k+1}{n}\right)\right|}{\sqrt{(k+1)/n}}\right\}\\
    & \le \sqrt{\frac{1}{n}} + 4\sqrt{2}\cdot \sqrt{\frac{1}{n}}\log\left(\frac{n}{\delta}\right)\le 8\sqrt{\frac{1}{n}}\log\left(\frac{n}{\delta}\right).
    \end{align*}
    Therefore, we have proved that with probability at least $1 - \frac{\delta}{2}$, 
    $$\frac{|F(t) - \hat{F}_n(t)|}{\sqrt{1/n\vee(t\wedge (1-t))}}\le 8\sqrt{\frac{1}{n}}\log\left(\frac{n}{\delta}\right)$$
    holds for every $0\le t\le \frac{1}{2}$. Similarly, we can prove that with probability at least $1 - \frac{\delta}{2}$, the above inequality holds for $\frac{1}{2}\le t\le 1$. Therefore, with probability at least $1 - \delta$, we have
    $$\sup_{0\le t\le 1}\frac{|F(t) - \hat{F}_n(t)|}{\sqrt{1/n\vee(t\wedge (1-t))}}\le 8\sqrt{\frac{1}{n}}\log\left(\frac{n}{\delta}\right).$$
    This completes the proof of this lemma.
   	
\end{proof}

\begin{proof}[Proof of Proposition \ref{prop-difference}]
    Suppose random variables $X\sim \mathbb{P}, Y\sim \mathcal{N}$ are independent. Then $X + Y\sim \mathbb{P}*\mathcal{N}$. We generate $n$ \emph{i.i.d.} samples $X_1, \cdots, X_n$; $Y_1, \cdots, Y_n$. Then $X_i + Y_i$ are $n$ \emph{i.i.d.} samples of $\mathbb{P}*\mathcal{N}$. We use $\hat{F}_n$ to denote the PDF of empirical measure $\hat{\mathbb{P}}_n = \frac{1}{n}\sum_{i=1}^n\delta_{X_i + Y_i}$. Then according to Lemma \ref{lem-CDF}, we have with probability $1 - \delta$, 
    $$\sup_{t\in\mathbb{R}} \frac{|F(t) - \hat{F}_n(t)|}{\sqrt{1/n\vee(F(t)\wedge (1-F(t)))}}\le \frac{8}{\sqrt{n}}\log\left(\frac{n}{\delta}\right).$$
    Hence Markov inequality indicates that 
    \begin{align*}
		\mathbf{P}\left(\exp\left(\sup_{t\in\mathbb{R}} \frac{\sqrt{n}}{16}\cdot\frac{|F(t) - \hat{F}_n(t)|}{\sqrt{1/n\vee(F(t)\wedge (1-F(t)))}} - \frac{\log (n)}{2}\right)\ge \frac{1}{\delta}\right)\le \delta^{2}. \numberthis \label{eq: inequality-p}
	\end{align*}
    Moreover, we notice that
    $$\mathbb{E}\left[\hat{F}_n(t)\Big|X_1, \cdots, X_n\right] = \mathbf{P}\left(\frac{1}{n}\sum_{i=1}^n (X_i + Y_i)\le t\Big|X_1, \cdots, X_n\right) = F_{n}(t),$$
    where $F_{n}$ is the CDF of $\mathbb{P}_n * \mathcal{N}$ with $\mathbb{P}_n = \frac{1}{n}\sum_{i=1}^n\delta_{X_i}$. We further notice $|F(t) - F_n(t)| = |F(t) - \mathbb{E}_{Y_{i}, 1\le i\le n}[\hat{F}_n(t)]|\le \mathbb{E}_{Y_{i}, 1\le i\le n}[F(t) - \hat{F}_n(t)]$.Hence according to Jensen's inequality,
	$$\exp\left(\frac{|F(t) - F_n(t)|}{\sqrt{1/n\wedge (F(t)\vee (1-F(t)))}}\right)\le \mathbb{E}\left[\exp\left(\sup_{t\in\mathbb{R}}\frac{|F(t) - \hat{F}_n(t)|}{\sqrt{1/n\wedge (F(t)\vee (1-F(t)))}}\right)\mid X_1, \cdots, X_n\right],$$
	which implies
    \begin{align*}
    	&\quad \mathbb{E}\left[\exp\left(\sup_{t\in\mathbb{R}} \frac{\sqrt{n}}{16}\cdot\frac{|F(t) - F_n(t)|}{\sqrt{1/n\vee(F(t)\wedge (1-F(t)))}} - \log(n)\right)\right]\\
		& \le \mathbb{E}\left[\mathbb{E}\left[\exp\left(\sup_{t\in\mathbb{R}} \frac{\sqrt{n}}{16}\cdot\frac{|F(t) - \hat{F}_n(t)|}{\sqrt{1/n\vee(F(t)\wedge (1-F(t)))}} - \frac{\log (n)}{2}\right)\right]\Bigg| X_1, \cdots, X_n\right]\\
		& = 1 + \int_1^\infty \mathbf{P}\left(\exp\left(\sup_{t\in\mathbb{R}} \frac{\sqrt{n}}{16}\cdot\frac{|F(t) - \hat{F}_n(t)|}{\sqrt{1/n\vee(F(t)\wedge (1-F(t)))}} - \frac{\log (n)}{2}\right)\ge r\right)dr\le 1 + \int_{1}^{\infty}\frac{1}{r^{2}}dr = 2,
    \end{align*}
    where in the last inequality we use \eqref{eq: inequality-p}. And according to Markov inequality, we have
    $$\mathbf{P}\left(\exp\left(\sup_{t\in\mathbb{R}} \frac{\sqrt{n}}{16}\cdot\frac{|F(t) - F_{n}(t)|}{\sqrt{1/n\vee(F(t)\wedge (1-F(t)))}} - \log (n)\right)\ge \frac{2}{\delta}\right)\le \delta.$$
    Therefore, with probability at least $1 - \delta$ we have
    $$\sup_{t\in\mathbb{R}}\frac{|F_{n}(t) - F(t)|}{\sqrt{1/n\vee(F(t)\wedge (1-F(t)))}}\le \frac{16}{\sqrt{n}}\log\left(\frac{2n}{\delta}\right).$$
\end{proof}

We are now ready to prove the upper bound part of Theorem \ref{thm-1D}
\begin{proof}[Proof of the Upper Bound in Theorem \ref{thm-1D}]
	\par In the following proof, we use $\mathcal{N}$ to denote the 1-dimensional standard normal distribution $\mathcal{N}$, and $T(\cdot)$ to denote the push-forward operator between $\PP * \mathcal{N}$ and $\PP_n * \mathcal{N}$ ($T(t) = F_n^{-1}(F(t))$, where $F, F_n$ are CDF of distribution $\PP * \mathcal{N}$ and $\PP_n * \mathcal{N}$ respectively). 
	\par First, as shown in \eqref{eq-conditioning} in the outline of the proof, we show that if $E$ is any event of probability at least $1 - \frac{C_E}{n^2}$ for some constant $C_E$ only depending on $C, K$, then we have
	\begin{equation}\label{eq-event}\mathbb{E}\left[W_2^2(\mathbb{P} * \mathcal{N}, \PP_n * \mathcal{N})\right]\le \mathbb{E}\left[W_2^2(\mathbb{P} * \mathcal{N}, \PP_n * \mathcal{N})|E\right] + \mathcal{O}\left(\frac{1}{n}\right).\end{equation}
	Actually we have
	$$\mathbb{E}\left[W_2^2(\mathbb{P} * \mathcal{N}, \PP_n * \mathcal{N})\right] = \mathbb{E}\left[W_2^2(\mathbb{P} * \mathcal{N}, \PP_n * \mathcal{N})\mathbf{1}_E\right] + \mathbb{E}\left[W_2^2(\mathbb{P} * \mathcal{N}, \PP_n * \mathcal{N})\mathbf{1}_{E^c}\right],$$
	and $ \mathbb{E}\left[W_2^2(\mathbb{P} * \mathcal{N}, \PP_n * \mathcal{N})\mathbf{1}_E\right] =  \mathbb{E}\left[W_2^2(\mathbb{P} * \mathcal{N}, \PP_n * \mathcal{N})|E\right]\mathbf{P}(E)\le \mathbb{E}\left[W_2^2(\mathbb{P} * \mathcal{N}, \PP_n * \mathcal{N})|E\right]$. As for the second term, according to Cauchy-Schwartz inequality we have
	\begin{align*}
		\mathbb{E}\left[W_2^2(\mathbb{P} * \mathcal{N}, \PP_n * \mathcal{N})\mathbf{1}_{E^c}\right]& \le \sqrt{\mathbb{E}[\mathbf{1}_E]\mathbb{E}[W_2^4(\mathbb{P} * \mathcal{N}, \PP_n * \mathcal{N})]} = \sqrt{\mathbf{P}(E^c)\mathbb{E}[W_2^4(\mathbb{P} * \mathcal{N}, \PP_n * \mathcal{N})]}\\
		& \le \frac{\sqrt{C_E}}{n}\sqrt{\mathbb{E}[W_2^4(\mathbb{P} * \mathcal{N}, \PP_n * \mathcal{N})]}.\end{align*}
	We further notice that according to the triangle inequality of W2 distance we have
	\begin{align*}\mathbb{E}[W_2^4(\mathbb{P} * \mathcal{N}, \PP_n * \mathcal{N})] \le \mathbb{E}[\left(W_2(\mathbb{P} * \mathcal{N}, \delta_0) + W_2(\delta_0, \mathbb{P}_n * \mathcal{N})\right)^4]\le \mathbb{E}[8W_2(\mathbb{P} * \mathcal{N}, \delta_0)^4 + 8W_2(\delta_0, \mathbb{P}_n * \mathcal{N})^4]\\
	= 8\mathbb{E}[(V_1+Z)^4] + 8\mathbb{E}[\mathbb{E}[(V_2+Z)^4|X_{1:n}]] = 64\mathbb{E}[V_1^4] + 64\mathbb{E}[\mathbb{E}[V_2^4|X_{1:n}]] + 128\mathbb{E}[Z^4] = \mathcal{O}(1),
	\end{align*}
	where we use $\delta_0$ to denote the delta distribution at 0, and $V_1\sim \PP, V_2\sim \PP_n, Z\sim \mathcal{N}$ are all independent. The last equation is because $\PP$ is $K$-subgaussian, hence all moments of $\PP$ are upper bounded by constant. Hence we have proved \eqref{eq-event}.
	
	\par Next, we notice that for any $n$ $i.i.d.$ samples $X_1, \cdots, X_n$, we have
	$$\mathbf{P}\left[\left\{|X_i|\le 2K\sqrt{2\log n}, \forall 1\le i\le n\right\}\right] = \left(1 - \mathbb{P}\left[|X_1|\ge 2K\sqrt{2\log n}\right]\right)^n\ge \left(1 - \frac{C}{n^4}\right)^n\ge 1 - \frac{C}{n^3}\ge 1 - \frac{C}{n^2},$$
	where we use the fact that $\mathbb{P}$ is a $K$-subgaussian distribution ($\mathbf{P}(|X|\ge t)\le C\exp(-t^2/(2K^2))$). Therefore, with probability at least $1 - \frac{C}{n^2}$ we have for $X'\sim\mathbb{P}_n$
	$$\mathbf{P}(|X'|\ge r))\le e\exp\left(-\frac{r^2}{2(2K\sqrt{\log n})^2}\right),$$
	which indicates that $\PP_n$ is $2K\sqrt{\log n}$ subgaussian with probability at least $1 - \frac{C}{n^2}$. We assume this event to be $E$, then $\mathbf{P}(E)\ge 1 - \frac{C}{n^2}$. In the following proof we all assume the event $E$, where we will deal with $E^c$ in the end. Noticing that for $X\sim \PP, Y\sim \PP_n$ we have
	$$\mathbf{P}(|X|\ge r)\le C\exp\left(-\frac{r^2}{2K^2}\right), \quad \mathbf{P}(|Y|\ge r)\le e\exp\left(-\frac{r^2}{2(2K\sqrt{\log n})^2}\right),$$
	according to Proposition \ref{prop-pq} we obtain that
	\begin{align*}
		|T(x) - x| & \lesssim |x| + 1 + \sqrt{\log n}(|x| + 1 + \sqrt{\log(|x| + 1)})\lesssim 1 + \sqrt{\log n}|x| \numberthis \label{eq-transport}
	\end{align*}
	for some positive constant $C_1, C_2$ only depending on $C, K$. We further notice that according to \eqref{eq-form-rho} we have
	$$\rho(x)\le \left(\frac{1}{\sqrt{2\pi}} + C\right)\exp\left(-\frac{x^2}{8K^2}\right).$$
	Hence when $|x|\ge 2K\sqrt{2\log n}$, we will have $\exp\left(-\frac{x^2}{8K^2}\right)\le \frac{1}{n}$ and hence
	$$\int_{|t| > 2K\sqrt{2\log n}}\rho(t)|T(t) - t|^2\lesssim \left(\frac{1}{\sqrt{2\pi}} + C\right)\int_{|t| > 2K\sqrt{2\log n}}\exp\left(-\frac{t^2}{8K^2}\right)(1 + \sqrt{\log n}|t|)^2dt = \tilde{\mathcal{O}}\left(\frac{1}{n}\right)$$

	\par Therefore, we only need to analyze the integral
	\begin{equation}\label{eq-integral}\int_{|t|\le 2K\sqrt{2\log n}}\rho(t)|T(t) - t|^2.\end{equation}
	In what follows, we will use the notation:
	\begin{align*}
		\overline{\rho}(t) & = \sup_{x\in [t-1, t+1]}\rho(x),\quad \underline{\rho}(t) = \inf_{x\in [t-1, t+1]}\rho(x)\qquad \text{and}\quad\Lambda(t) = \sup_{x\in [t-1, t+1]}|F(x) - F_n(x)|.
	\end{align*}
	The key idea to bound the integral in \eqref{eq-integral} is the following observation from Proposition \ref{prop-cdf-pdf}: if $\Lambda(t)\le \underline{\rho}(t)$, then we have $|T(t) - t|\le \frac{\Lambda(t)}{\underline{\rho}(t)}$, which implies that
	$$\rho(t)|T(t) - t|^2\le \frac{\Lambda(t)^2}{\rho(t)}\cdot \left(\frac{\rho(T)}{\underline{\rho}(t)}\right)^2.$$
	In the following, we use the concentration proposition (Proposition \ref{prop-difference}) to divide the interval $[-2K\sqrt{\log n}, 2K\sqrt{\log n}]$ into the set where $\Lambda(t)\le \underline{\rho}(t)$ where the integral can be bounded from the above inequality, and the set where $\rho(t)$ is very small hence $\rho(t)|T(t) - t|^2$ won't have much effect in the integral.
	
	\par According to Proposition \ref{prop-difference}, with probability at least $1 - \frac{1}{n^2}$ we have
	$$\sup_{t\in\mathbb{R}}\frac{|F(t) - F_{n}(t)|}{\sqrt{1/n\vee \min\{F(t), 1-F(t)\}}}\le \frac{16}{\sqrt{n}}\log\left(2n^3\right).$$
	We assume this event to be $E_1$, where $\mathbf{P}(E_1)\ge 1 - \frac{1}{n^2}$. In the rest of the proof we assume $E_1$ and will deal with $E_1^c$ in the end. Then we have
	\begin{align*}
		\Lambda(t) = \sup_{x\in [t-1, t+1]}|F(x) - F_n(x)| & \le \frac{16\log (2n^3)}{\sqrt{n}}\sup_{x\in [t-1, t+1]}\sqrt{\frac{1}{n}\vee\min\{F(x), 1-F(x)\}}\\
		& = \frac{16\log (2n^3)}{\sqrt{n}}\sqrt{\frac{1}{n}\vee\sup_{x\in [t-1, t+1]}\min\{F(x), 1-F(x)\}}.
	\end{align*}
	According to Proposition \ref{prop-PDF-CDF}, for any $0 < \epsilon < \beta$, $\exists M = M(K, C, \epsilon)\ge 1$ such that $\min\{F(x), 1-F(x)\}\le M\rho(r)^{\beta - \epsilon},$ which indicates that
	$$\Lambda(t)\le \frac{16\log (2n^3)}{\sqrt{n}}\sqrt{\frac{1}{n}\vee \sup_{x\in [t-1, t+1]}M\rho(x)^{\beta - \epsilon}} = \frac{16\log (2n^3)}{n}\vee\frac{16\sqrt{M}\log (2n^3)}{\sqrt{n}}\overline{\rho}(t)^{\frac{\beta - \epsilon}{2}}$$
	 Next we will upper bound $\frac{\overline{\rho}(t)}{\rho(t)}$ and also $\frac{\rho(t)}{\underline{\rho}(t)}$ from the following observation: Noticing that for $S\sim \PP$ we have 
	 $$\EE[S] = \int_{-\infty}^{\infty}x\eta(x)dx\le \int_{-\infty}^\infty |x|\eta(x)dx = \int_0^\infty \mathbf{P}(S\ge r)dr\le \int_0^\infty C\exp\left(-\frac{r^2}{2K^2}\right)dr = \frac{CK\sqrt{2\pi}}{2}\le 2CK,$$
	hence according to~\cite[Prop. 2]{polyanskiy2016wasserstein}, we obtain that $\PP * \mathcal{N}$ is $(3, 8CK)$-regular, which indicates that for $|t|\le 2K\sqrt{2\log n}$ and $\forall x\in [t-1, t+1]$, 
	 \begin{align*}
	 	& \frac{\rho(x)}{\rho(t)}\le \exp\left(3(|t|+1) + 8CK\right)\le \exp\left(6K\sqrt{2\log n} + 3 + 8CK\right)\triangleq L(n)\\
		& \frac{\rho(x)}{\rho(t)}\ge \exp\left(-3(|t|+1) - 8CK\right)\le \exp\left(-6K\sqrt{2\log n} - 3 - 8CK\right) = \frac{1}{L(n)}.
	\end{align*}
	Hence we have $1\le \frac{\overline{\rho}(t)}{\rho(t)}, \frac{\rho(t)}{\underline{\rho}(t)}\le L(n)$. Therefore, when 
	$$\rho(t)\ge \frac{16\log (2n^3)L(n)}{n}\vee \left(\frac{256M\log^2(2n^3)}{n}L(n)^{2+\beta-\epsilon}\right)^\frac{1}{2 - \beta + \epsilon}\triangleq Q(n),$$
	we will have
	$$\Lambda(t)\le \frac{16\log (2n^3)}{n}\vee\frac{16\sqrt{M}\log (2n^3)}{\sqrt{n}}\overline{\rho}(t)^{\frac{\beta - \epsilon}{2}}\le \underline{\rho}(t),$$
	which, according to Proposition \ref{prop-cdf-pdf} with $h = 1$, we have $|T(t) - t|\le \frac{\Lambda(t)}{\underline{\rho}(t)}\le L(n)\frac{\Lambda(t)}{\rho(t)}$. Therefore, noticing that $\beta\le 1$, we have
	\begin{align*}
		&\quad \int_{\{t|\rho(t)\ge Q(n), |t|\le 2K\sqrt{\log n}\}}\rho(t)|T(t) - t|^2dt\le L(n)^2\int_{\{t|\rho(t)\ge Q(n), |t|\le 2K\sqrt{\log n}\}}\frac{\Lambda(t)^2}{\rho(t)}dt\\
		& \le 4KL(n)^2\sqrt{\log n}\cdot \max_{\rho(t)\ge Q(n)}\frac{\Lambda(t)^2}{\rho(t)}\le 4KL(n)^2\sqrt{\log n}\cdot\max_{\rho(t)\ge Q(n)}\left\{\frac{256\log^2 (2n^3)}{n^2\rho(t)}\vee\frac{256M\log^2(2n^3)}{n\rho(t)^{1+\epsilon-\beta}}\right\}\\
		& \le 4KL(n)^2\sqrt{\log n}\cdot\left(\frac{256\log^2 (2n^3)}{n^2Q(n)}\vee\frac{256M\log^2(2n^3)}{nQ(n)^{1+\epsilon-\beta}}\right)\\
		& \le 4KL(n)^2\sqrt{\log n}\cdot\left(\frac{16\log (2n^3)}{nL(n)}\vee\left(\frac{256M\log^2(2n^3)}{n}\right)^\frac{1}{2-\beta+\epsilon}L(n)^{-\frac{(2+\beta-\epsilon)(1+\epsilon-\beta)}{2-\beta+\epsilon}}\right).
	\end{align*}
	Further noticing that for any $\epsilon_1 > 0$, we have $L(n) = \mathcal{O}\left(n^{\epsilon_1}\right)$. Hence for any $\epsilon' > 0$, we have
	$$\int_{\{t|\rho(t)\ge Q(n), |t|\le 2K\sqrt{\log n}\}}\rho(t)|T(t) - t|^2dt = \mathcal{O}\left(n^{-\frac{1}{2 - \beta + \epsilon} + \epsilon'}\right).$$
	
	\par As for those $t$ with $\rho(t) < Q(n)$, according to \eqref{eq-transport} we have estimation
	\begin{align*}
		&\quad \int_{\{t|\rho(t) < Q(n), |t|\le 2K\sqrt{\log n}\}}\rho(t)|T(t) - t|^2dt\le \int_{\{t|\rho(t) < Q(n), |t|\le 2K\sqrt{\log n}\}}\rho(t)(C_1 + C_2\sqrt{\log n}|t|)^2dx\\
		&\le 4K\sqrt{\log n}\cdot Q(n)(C_1 + 2KC_2\log n)^2 = Q(n)\cdot \tilde{\mathcal{O}}(1) = \mathcal{O}\left(n^{-\frac{1}{2 - \beta + \epsilon} + \epsilon'}\right)
	\end{align*}
	Combine these two estimation together, we obtain that assuming event $E, E_1$, for any $\epsilon, \epsilon' > 0$, we have $\int_{|t|\le 2K\sqrt{n}}\rho(t)|T(t) - t|^2dt = \mathcal{O}\left(n^{-\frac{1}{2 - \beta + \epsilon} + \epsilon'}\right)$ and hence
	$$\mathbb{E}\left[W_2^2(\mathbb{P} * \mathcal{N}, \PP_n * \mathcal{N})|E\cap E_1\right] = \int_{-\infty}^{\infty}\rho(t)|T(t) - t|^2dt = \mathcal{O}\left(n^{-\frac{1}{2 - \beta + \epsilon} + \epsilon'}\right) + \tilde{\mathcal{O}}(\frac{1}{n}) = \mathcal{O}\left(n^{-\frac{1}{2 - \beta + \epsilon} + \epsilon'}\right).$$
	
	\par Finally we notice that $\mathbf{P}(E^c\cup E_1^c) = \frac{C + 1}{n^2}$, according to \eqref{eq-event} we have
	$$\mathbb{E}\left[W_2^2(\mathbb{P} * \mathcal{N}, \PP_n * \mathcal{N})\right]\le \mathbb{E}\left[W_2^2(\mathbb{P} * \mathcal{N}, \PP_n * \mathcal{N})|E\cap E_1\right] + \mathcal{O}\left(\frac{1}{n}\right) = \mathcal{O}\left(n^{-\frac{1}{2 - \beta + \epsilon} + \epsilon'}\right).$$
	Since $\epsilon$ and $\epsilon'$ can be chosen to be arbitrary small positive number, and $2\alpha = \frac{(1 + K^2)^2}{2(1 + K^4)} = \frac{1}{2-\beta}$, we have for any $\epsilon > 0$, 
	$$\mathbb{E}\left[W_2^2(\mathbb{P} * \mathcal{N}, \PP_n * \mathcal{N})\right] = \mathcal{O}\left(n^{-2\alpha + \epsilon}\right).$$ 	
\end{proof}

\section{Proof of Theorem \ref{thm-KL}}\label{thm-KL-part}
To prove Theorem \ref{thm-KL}, we will leverage the concepts of R\'{e}nyi divergence and R\'{e}nyi mutual information, defined as follows~\cite[Section 7.12]{itbook}:
\begin{definition}[R\'{e}nyi Divergence and R\'{e}nyi Mutual Information~\cite{renyi1961measures}]\label{def-renyi}
	Assume random variables $(X, Y)$ have joint distribution $P_{X, Y}$. For any $\lambda > 1$, the R\'{e}nyi divergence between two distributions $\PP$ and $\QQ$ is defined as 
    $$D_{\lambda}(\PP\|\QQ) \triangleq \frac{1}{\lambda - 1}\log \mathbb{E}_{\QQ}\left[\left(\frac{d\PP}{d\QQ}\right)^{\lambda}\right].$$
    And the R\'{e}nyi Mutual Information of order $\lambda$ are defined as
	$$I_\lambda(X; Y) \eqdef D_\lambda(P_{X,Y} \| P_X \otimes P_Y),$$
	where we use $P_X, P_Y$ to denote the marginal distribution with respect to $X$ and $Y$, and $P_X\otimes P_Y$ denotes the joint distribution of $(X', Y')$ where $X'\sim P_X, Y'\sim P_Y$ are independent to each other. 
\end{definition}

A key result we need in order to establish Theorem~\ref{thm-KL} is the following lemma which upper bounds smoothed KL divergence by R\'{e}nyi mutual information.
\begin{lemma}\label{lem-KL-I-lambda}
	Suppose $(X, Y)\sim \PP_{X, Y}$, with marginal distributions $\PP_X, \PP_Y$.  Let $\PP_{n}$ be an empirical version of $\PP_{X}$ generated with $n$ samples. Then for every $1 < \lambda \le 2$, we have
	\begin{equation}\label{eq:softcov_kl}
		\EE[D_{KL}(\PP_{Y|X} \circ \PP_n \| \PP_Y)] \le {1\over \lambda -1} \log( 1 + \exp\{(\lambda -1)(I_\lambda(X;Y) - \log n)\})\,.
	\end{equation}	
\end{lemma}

\begin{proof} 
\par According to \cite{van2014renyi}, for any distribution $\PP, \QQ$, the function $D_{\lambda}(\PP\|\QQ)$ with respect to $\lambda\in (1, 2]$ is non-decreasing, where $D_{\lambda}$ is the R\'{e}nyi divergence defined in Definition \ref{def-renyi}. Hence noticing from \cite{van2014renyi} that for any distribution $\PP, \QQ$, $\lim_{\lambda\to 1}D_{\lambda}(\PP\|\QQ) = D_{KL}(\PP\|\QQ)$, we have
$$D_{KL}(\PP_{Y|X} \circ \PP_n \| \PP_Y)\le D_{\lambda}(\PP_{Y|X} \circ \PP_n \| \PP_Y).$$
Therefore, it is sufficient to prove that for any $1 < \lambda \le 2$,  
$$\EE[D_{\lambda}(\PP_{Y|X} \circ \PP_n \| \PP_Y)] \le {1\over \lambda -1} \log( 1 + \exp\{(\lambda -1)(I_\lambda(X;Y) - \log n)\}).$$ 
We suppose the $n$ samples obtained in $\PP_{n}$ to be $X_{1}, \cdots, X_{n}$, which satisfies that $(X_1, \cdots, X_n)\indep Y$. According to the definition of R\'{e}nyi divergence, R\'{e}nyi mutual information and also the Jensen's inequality, we see that
\begin{align} \EE[D_\lambda(\PP_{Y|X} \circ \PP_n \| \PP_Y)] & =  {1\over \lambda -1} \EE\left[\log
\EE\left[ \left\{ {d(\PP_{Y|X} \circ \PP_n)(Y) \over d\PP_Y(Y)} \right\}^\lambda \right]\Bigg|X_{1:n}\right]
\label{eq:sco_3}\\
		& \le {1\over \lambda -1} \log
\EE\left[ \left( {d(\PP_{Y|X} \circ \PP_n)(Y) \over d\PP_Y(Y)} \right)^\lambda \right].\nonumber
\end{align}
Then we introduced the channel $\PP_{\bar Y|X_{1:n}} = {1\over n}
\sum_{i=1}^n \PP_{Y|X=X_i}$ and we let $\PP_{X_{1:n}, \bar{Y}} = \PP_{\bar{Y}|X_{1:n}}\circ \PP_{X_{1:n}}$, where $\PP_{X_{1:n}} = \PP_X^{\otimes n}$ is the probability law of $X_{1:n}$. We notice that the marginal distribution of $\PP_{X_{1:n}, \bar{Y}}$ with respect to $\bar{Y}$ is exactly $\PP_Y$. If we let $(X_{1:n}, \bar{Y})\sim \PP_{X_{1:n}}\otimes \PP_{Y}$, then we obtain that
\begin{align*}
	I_\lambda(X_{1:n}; \bar Y) & = \frac{1}{\lambda - 1}\log\mathbb{E}\left[\left(\frac{d\PP_{X_{1:n}, \bar{Y}}(X_{1:n}, Y)}{d\left[\PP_{X_{1:n}}\otimes \PP_{Y}(X_{1:n}, Y)\right]}\right)^{\lambda}\right] = {1\over \lambda -1} \log \EE\left[ \left\{ {d\PP_{Y|X_{1:n}}(Y|X_{1:n}) \over d\PP_Y(Y)} \right\}^\lambda \right]\\
	& = {1\over \lambda -1} \log \EE\left[\EE\left[ \left\{ {d(\PP_{Y|X} \circ  \PP_n)(Y) \over d\PP_Y(Y)} \right\}^\lambda\Bigg| X_{1:n}\right]\right] = {1\over \lambda -1} \log\EE\left[ \left( {d(\PP_{Y|X} \circ \PP_n)(Y) \over d\PP_Y(Y)} \right)^\lambda \right]\\
	& \ge \EE[D_\lambda(\PP_{Y|X} \circ \PP_n \|\PP_Y)].
\end{align*}
Hence we only need to analyze $I_\lambda(X_{1:n}; \bar Y)$. And we need to upper bound
\begin{equation}\label{eq:sco_2}
	\EE\left[ \left\{ {d\PP_{Y|X_{1:n}}(Y|X_{1:n}) \over d\PP_Y(Y)} \right\}^\lambda \right] = \EE\left[ \left\{ {1\over n} \sum_{i=1}^n {d\PP_{Y|X}(Y|X_i)\over d\PP_Y(Y)}\right\}^{\lambda} \right]\,.
\end{equation} 
Moreover, noticing that $(a+b)^{\lambda -1} \le a^{\lambda -1} +b^{\lambda -1}$ holds for $a, b > 0$ and $1 < \lambda\le 2$, we have that for any $n$ $i.i.d.$ non-negative random variables $B_i$ ($1\le i\le n$), 
\begin{align*}
	\EE \left[B_i \left(B_i + \sum_{j\neq i} B_j\right)^{\lambda -1}\right] & \le  \EE [B_{i}\cdot B_{i}^{\lambda - 1}] +  \EE\left[B_{i}\cdot\left(\sum_{j\neq i} B_{j}\right)^{\lambda-1}\right] = \EE [B_{1}^{\lambda}] + \mathbb{E}[B_{i}]\cdot \EE\left[\left(\sum_{j\neq i} B_{j}\right)^{\lambda-1}\right]\\
	& \le \EE [B_{1}^{\lambda}] + \mathbb{E}[B_{1}]\cdot \left(\sum_{j\neq i}\mathbb{E}[B_{j}]\right)^{\lambda-1}\ = \EE [B_{1}^{\lambda}] + \mathbb{E}[B_{1}]\cdot \left((n-1)\mathbb{E}[B_{1}]\right)^{\lambda-1},
\end{align*}
where in the second inequality we use the Jensen's inequality. Therefore, summing up the above inequality for $1\le i\le n$, we have
$$\EE\left[ \left\{\sum_{i=1}^n B_i\right\}^\lambda \right] \le n \EE[B_{1}^\lambda] + n\cdot (n-1)^{\lambda-1}\left(\EE[B_{1}]\right)^\lambda\le n \EE[B_{1}^\lambda] + n^{\lambda}\left(\EE[B_{1}]\right)^\lambda\,.$$
This puts us into a well-known setting of Rosenthal-type inequalities, which is known to be essentially tight~\cite{schechtman2011extremal}.

\par Next, since $Y\indep (X_1, \cdots, X_n)$, for every fixed $Y$, random variables $\frac{d\PP_{Y|X}(Y|X_i)}{d\PP_{Y}(Y)}$ are $i.i.d$. Hence choosing $B_{i} = \frac{d\PP_{Y|X}(Y|X_i)}{d\PP_{Y}(Y)}$, we obtain that
2\begin{align*}
	&\quad \EE\left[ \left\{ {1\over n} \sum_i {d\PP_{Y|X}(Y|X_i)\over d\PP_Y(Y)}\right\}^{\lambda}\Bigg|Y\right]\le n^{-\lambda}\cdot \EE\left[ \left\{\sum_i {d\PP_{Y|X}(Y|X_i)\over d\PP_Y(Y)}\right\}^{\lambda}\Bigg|Y\right]\\
	&\le n^{-\lambda}\cdot \left(n\cdot \EE\left[ \left\{  {d\PP_{Y|X}(Y|X)\over d\PP_Y(Y)}\right\}^{\lambda}\bigg|Y\right] + n^{\lambda}\cdot \left(\EE\left[{d\PP_{Y|X}(Y|X)\over d\PP_Y(Y)}\bigg|Y\right]\right)^\lambda\right)\\
	&\le n^{1-\lambda} \EE\left[ \left\{  {d\PP_{Y|X}(Y|X)\over d\PP_Y(Y)}\right\}^{\lambda}\bigg| Y\right] + \left(\EE\left[{d\PP_{Y|X}(Y|X)\over d\PP_Y(Y)}\bigg|Y\right]\right)^\lambda.
\end{align*}
Using the fact that $X\indep Y$ and hence $\mathbb{E}[\PP_{Y|X}(Y|X)|Y] = \int_X P_{Y|X}(Y|X)d\PP_X(X) = \int_X d\PP_{X, Y}(X, Y) = \PP_Y(Y)$, we notice that for any given $Y$,
$$\EE\left[{d\PP_{Y|X}(Y|X)\over d\PP_Y(Y)}\bigg|Y\right] = \frac{d\mathbb{E}[\PP_{Y|X}(Y|X)]}{d\PP_Y(Y)}\bigg|_Y = \frac{d\PP_Y(Y)}{d\PP_Y(Y)}\bigg|_Y = 1.$$
Therefore, we can upper bound \eqref{eq:sco_2} as
\begin{align*}
	&\quad \EE\left[ \left\{ {1\over n} \sum_{i=1}^n {d\PP_{Y|X}(Y|X_i)\over d\PP_Y(Y)}\right\}^{\lambda} \right] = \EE\left[\EE\left[ \left\{ {1\over n} \sum_{i=1}^n {d\PP_{Y|X}(Y|X_i)\over d\PP_Y(Y)}\right\}^{\lambda} \right]\Bigg| Y\right]\\
	&\le n^{1-\lambda} \EE\left[\EE\left[ \left\{  {d\PP_{Y|X}(Y|X)\over d\PP_Y(Y)}\right\}^{\lambda}\bigg| Y\right]\bigg|Y\right] + \EE\left[\left(\EE\left[{d\PP_{Y|X}(Y|X)\over d\PP_Y(Y)}\bigg|Y\right]\right)^\lambda\bigg| Y \right]\\
	&\le n^{1-\lambda} \EE\left[ \left\{  {d\PP_{Y|X}(Y|X)\over d\PP_Y(Y)}\right\}^{\lambda} \right] + 1 = n^{1 - \lambda}\cdot \exp\left((\lambda - 1)I_{\lambda}(X; Y)\right) + 1.
\end{align*}
This implies that
$$ I_\lambda(X_{1:n}; \bar Y) \le {1\over \lambda -1} \log \left( 1+ n^{1-\lambda} \exp\{(\lambda -1)
I_\lambda(X;Y)\} \right)\,,$$
which together with~\eqref{eq:sco_3} recovers~\eqref{eq:softcov_kl}.

\end{proof}

\begin{remark} Hayashi~\cite{hayashi2006general} upper bounds the LHS of~\eqref{eq:softcov_kl}
with 
$$ {\lambda \over \lambda - 1} \log \left(1 + \exp\left\{{\lambda -1\over \lambda} (K_\lambda(X;Y) - \log
n)\right\}\right)\,,$$
where $K_\lambda(X;Y) = \inf_{\QQ_Y} D_\lambda(\PP_{X,Y} \| \PP_X \otimes\QQ_Y)$ is the so-called Sibson-Csiszar
information, cf.~\cite{sibson1969information}. This bound, however, does not have the right rate of convergence
as $n\to\infty$, at least for $\lambda=2$ as comparison with Prop. 5 in~\cite{goldfeld2020convergence}.
We note that~\cite{hayashi2006general,han1993approximation} also contain bounds on $\mathbb{E}[\mathrm{TV}(\PP_{Y|X} \circ \PP_n, \PP_Y)]  $
which do not assume existence of $\lambda>1$ moment of ${\PP_{Y|X}\over \PP_Y}$ and instead rely on
the distribution of $\log {d\PP_{Y|X}\over d\PP_Y}$. 

\end{remark}

On top of the previous lemma, we use the following lemma to provide an upper bound to the R\'{e}nyi mutual information.
\begin{lemma}\label{lem-I-lambda}
	Suppose random variables $X\sim
\mathbb{P}, Z\sim \mathcal{N}(0, \sigma^{2}I_{d})$ are independent of each other. Let $Y = X +
Z$. Fix $1 < \lambda < 2$ and let $l=l(\lambda) = {\lambda -1 \over 2-\lambda}(d+1)$. If the random variable $X\sim\PP$ has finite $l$-th moment $\mathcal{M}
\eqdef \EE[\|X\|^l]^{2-\lambda}$, then for any $\sigma > 0$ :
$$ I_\lambda(X; Y) \le  \frac{1}{\lambda - 1}\log\left(C\mathcal{M}\right)\,,$$
where $C = C(\sigma)>0$ only depends on $\sigma$.
Moreover, if $\PP$ is a $K$-subgaussian distribution, we have for all $1<\lambda < 2$  
$$I_{\lambda}(X; Y)\le \frac{1}{\lambda - 1}\log\left(\frac{C'}{(2-\lambda)^{d}}\right)$$
for some constant $C' = C(K, \sigma) > 0$. 
\end{lemma}

\begin{proof}
	According to the definition of R\'{e}nyi divergence, we have 
	$$I_\lambda(X; Y) = \frac{1}{\lambda - 1}\log \left(C_0\EE\left[\int{\rho_{Y|X}^{\lambda}(\mathbf{y}|X)\over \rho_Y(\mathbf{y})^{\lambda-1}}d\mathbf{y}\right]\right),$$
	for some positive constant $C_0$, where $\rho_{Y|X}(\my|X) = \exp(-{\|\my-X\|^2\over 2\sigma^2})$ and $\rho(\my) = \mathbb{E}\left[\exp\left(-\frac{\|\my - X\|_{2}^{2}}{2\sigma^{2}}\right)\right].$
	Therefore, we only need to prove $\EE[\int{\rho_{Y|X}^{\lambda}(\mathbf{y}|X)\over \rho_Y(\mathbf{y})^{\lambda-1}}d\mathbf{y}] \lesssim  \mathcal{M}$
	for distributions $\PP$ with finite $l$-th moment, 
	and $\EE[\int{\rho_{Y|X}^{\lambda}(\mathbf{y}|X)\over \rho_Y(\mathbf{y})^{\lambda-1}}d\mathbf{y}] \lesssim  \frac{1}{(2-\lambda)^{d}}$ for $K$-subgaussian distribution $\PP$.
	\par We write $\mathbb{R} = \bigcup_i c_i$ as a union of cubes of diameter 2. For any $s\in c_i$, according to \eqref{eq:brt_2} we have
	$$\rho_Y(\mathbf{y})\ge \exp\left(-\frac{6}{\sigma^{2}}\right) \mathbf{P}(X\in c_i)\exp\left(-{3\|\my - s\|^{2}\over 4\sigma^{2}}\right),$$
	which indicates that
	\begin{align*}
		\frac{\exp\left(-\lambda\|\my - s\|_{2}^{2}/(2\sigma^{2})\right)}{\rho_Y(\my)^{\lambda - 1}} & \le \exp\left(\frac{6(\lambda - 1)}{\sigma^{2}}\right)\mathbf{P}(X\in c_i)^{1-\lambda} \exp\left(-{(3 - \lambda)\|\my -s \|^2\over 4\sigma^2}\right)\\
		& \le \exp\left(\frac{6}{\sigma^{2}}\right)\mathbf{P}(X\in c_i)^{1-\lambda} \exp\left(-{\|\my -s \|^2\over 4\sigma^2}\right)
	\end{align*}
	after noticing the fact that $1\le \lambda < 2$. Therefore for any $s\in\mreals^d$ we have
	$$\int_{\mreals^d}  \frac{\exp\left(-\lambda\|\my -
		s\|_{2}^{2}/(2\sigma^{2})\right)}{\rho_Y(\my)^{\lambda - 1}}d\my\lesssim \mathbf{P}(X\in c_i)^{1-\lambda}\int_{\mreals^d}\exp\left(-\frac{\|\mathbf{u}\|^2}{4\sigma^2}\right)d\mathbf{u}\lesssim \mathbf{P}(X\in c_i)^{1-\lambda},$$
	where we use $\lesssim$ to hide constant factors depending on $\sigma, d$. Taking the expectation over $X$, we obtain that
	\begin{equation}\label{eq-2-lambda}\EE\left[ {\rho_{Y|X}^{\lambda}(\my|X)\over \rho_Y(\my)^{\lambda-1}}\right]\lesssim \sum_i \mathbf{P}(X\in c_i)^{2-\lambda}.\end{equation}
	
	\par Next, we use $L_r$ 
	to denote the set of cubes whose centers belong to $\{r-1\le \|x_i\| < r\}$. Then we have $|L_r| = \mathcal{O}(r^{d-1})$. We further let $p_r = \sum_{c_i\in M_r} \mathbf{P}(X\in c_i)$, then according to Jensen's inequality we obtain that
	$$\sum_{c_i\in L_r}\mathbf{P}(X\in c_i)^{2-\lambda}\le |L_r|\cdot \left(\frac{1}{|L_r|}\sum_{c_i\in L_r}\mathbf{P}(X\in c_i)\right)^{2-\lambda} = |L_r|\cdot \left(\frac{p_r}{|L_r|}\right)^{2-\lambda} = |L_r|^{\lambda - 1}p_r^{2-\lambda}.$$
	Next recall that $l = \frac{(\lambda - 1)(d + 1)}{2 - \lambda}$. Assuming $\mathcal{M}_{l} < \infty$, we have for any $1 < \lambda < 2$,
	\begin{align*}
		\sum_{i=1} \mathbf{P}(X\in c_i)^{2-\lambda} & = \sum_{r=1}^{\infty}\sum_{c_i\in L_r}\mathbf{P}(X\in c_i)\lesssim \sum_{r=1}^\infty r^{(\lambda - 1)(d-1)}p_r^{2-\lambda}\\
		& \le \left(\sum_{r=1}^{\infty} r^\frac{(\lambda - 1)(d-1) + 2(\lambda - 1)}{2-\lambda}p_r\right)^{2-\lambda}\left(\sum_{r=1}^\infty\frac{1}{r^2}\right)^{\lambda - 1}\lesssim 1 + \EE[\|X\|^l]^{2-\lambda} = 1 + \mathcal{M}_l
	\end{align*}
	where in the second last inequality we use the H\"{o}lder inequality. Specifically for $K$-subgaussian cases, we notice that $p_{r}\lesssim\exp\left(-\frac{r^{2}}{2K^{2}}\right).$
Therefore, we obtain that
\begin{align*}
	&\quad \sum_{i=1} \mathbf{P}(X\in c_i)^{2-\lambda} = \sum_{r=1}^{\infty}\sum_{c_{i}\in L_{r}}\mathbf{P}(X\in c_i)^{2-\lambda}\le \sum_{r=1}^{\infty}|L_{r}|p_{r}^{2-\lambda}\lesssim \sum_{r=1}^{\infty}r^{d-1}\exp\left(-\frac{(2-\lambda)r^{2}}{2K^{2}}\right)\\
	& \le \sum_{r=1}^{\infty}r^{d-1}\exp\left(-\frac{(2-\lambda)r}{2K^{2}}\right)\lesssim\left(1 - \exp\left(-\frac{2-\lambda}{2K^{2}}\right)\right)^{-d}\le \frac{(2K^2)^d}{(2-\lambda)^{d}}\lesssim  \frac{1}{(2-\lambda)^{d}},
\end{align*}
where we use $\lesssim$ to hide coefficients depending on $K, d$.
	\par Based on these two upper bounds on $\sum \mathbf{P}(X\in c_{i})^{2-\lambda}$, \eqref{eq-2-lambda} yields the desired bounds on $I_{\lambda}(X; Y)$.
\end{proof}

\par We are now ready to prove Theorem \ref{thm-KL}.

\begin{proof}[Proof of Theorem \ref{thm-KL}]
	We consider $X\sim \mathbb{P}, Z\sim\mathcal{N}(0, \sigma^{2}I_{d}), X\indep Z$ and $Y = X + Z$. Then conditioned on $X$, we have $Y\sim \mathcal{N}(X, \sigma^{2}I_{d})$, which indicates that $\PP_{Y|X}\circ \mathbb{P}_{n}\sim \mathbb{P}_{n} * \mathcal{N}(0, \sigma^{2}I_{d})$. Therefore, adopting Lemma \ref{lem-KL-I-lambda} and Lemma \ref{lem-I-lambda}, we obtain that for any $1 < \lambda < 2$,
	\begin{align*}
		& \hspace{-0.5cm}\mathbb{E}[D_{KL}(\mathbb{P}_{n} * \mathcal{N}(0, \sigma^{2}I_d)\|\mathbb{P} * \mathcal{N}(0, \sigma^{2}I_d))]\le \frac{1}{\lambda - 1}\log (1 + \exp((\lambda - 1)(I_{\lambda}(X; Y) - \log n)))\\
		& \le \frac{1}{\lambda - 1}\cdot \exp((\lambda - 1)(I_{\lambda}(X; Y) - \log n))\le \frac{C}{(\lambda - 1)n^{\lambda - 1}(2 - \lambda)^{d}}.
	\end{align*}
	Choosing $\lambda = 2 - \frac{1}{\log n}$, and noticing that 
	$$n^{\lambda - 1} = n^{-\frac{1}{\log n} + 1} = x\cdot \exp\left(-\log n\cdot \frac{1}{\log n}\right) = \frac{n}{e},$$
	we have
	$$\mathbb{E}[D_{KL}(\mathbb{P}_{n} * \mathcal{N}(0, \sigma^{2}I_d)\|\mathbb{P} * \mathcal{N}(0, \sigma^{2}I_d))]\le \frac{Ce(\log n)^{d}}{(1 - 1/\log n)n} = \mathcal{O}\left(\frac{(\log n)^{d}}{n}\right).$$
	Hence \eqref{eq-upper-bound1} holds, which implies the upper bound part of Theorem \ref{thm-1D}.
\end{proof}

\section*{Acknowledgments}
 This work was supported by the MIT-IBM Watson AI Lab.
AB acknowledges support from the National Science Foundation Graduate Research Fellowship under
Grant No. 1122374.
AR acknowledges support from the Simons Foundation and the NSF through awards DMS-2031883 and DMS-1953181.

\bibliographystyle{alpha}
\bibliography{References}

\newpage
\appendix

\section{$W_2$-rates and log-Sobolev inequalities for Gaussian-smoothed Bernoulli measures}\label{sec: app-examples}

\subsection{Convergence Rate of  Wasserstein Distance}\label{sec-bernoulli}
 Our proof will rely on the following auxiliary lemma.
\begin{lemma}\label{lem2}
    Suppose two 1-dimensional distribution $\mu, \nu$ with CDFs $F_\mu, F_\nu$ satisfy $F_\mu(t)\ge F_\nu(t+2)$, then we have
    $$W_2(\mu, \nu)^{2}\ge \mathbf{P}(Y\in [t+1, t+2]), \quad Y\sim \nu.$$
\end{lemma}
\begin{proof}
    According to the definition of $W_2$ distance, we have $W_2(\mu, \nu)^2 = \inf_{\gamma}\int (x-y)^2d\gamma(x, y)$, where the infimum is over all possible coupling between $\mu$ and $\nu$. Hence we only need to prove that for any coupling $\gamma$ between $\mu$ and $\nu$, we always have $\mathbf{P}(|X - Y|\ge 1)\ge \mathbf{P}(Y\in [t+1, t+2])$, where $(X, Y)\sim \gamma$. Notice that $F_\mu(t)\ge F_\nu(t+2)$, we have $\mathbf{P}(X\le t)\ge \mathbf{P}(Y\le t+2) = \mathbf{P}(Y < t+1) + \mathbf{P}(Y\in [t+1, t+2])$, which indicates that
    \begin{align*}
        \mathbf{P}(|X - Y|\ge 1) & \ge \mathbf{P}(X\le t, Y\ge t+1)\ge \mathbf{P}(X\le t) + \mathbf{P}(Y\ge t+1) - 1\\
        & \ge \mathbf{P}(Y < t+1) + \mathbf{P}(Y\in [t+1, t+2]) + \mathbf{P}(Y\ge t+1) - 1 = \mathbf{P}(Y\in [t+1, t+2]).
    \end{align*}
\end{proof}

\begin{proof}[Proof of Proposition \ref{prop-bernoulli}]

    Given $h > 0$, we assume $\mathbb{P}_{h, n} = (1 - \hat{p}_{h})\delta_0 + \hat{p}_{h}\delta_h$, where $\hat{p}_{h} = \frac{1}{n}\left(\sum_{k=1}^n \mathbf{1}_{X_k = h}\right)$, and $X_1, \cdots, X_n\sim \mathbb{P}_{h}$ are i.i.d. In the following proof, when there is no danger of confusion, we use $\tilde{F}_{n, \sigma}, F_\sigma, \Phi_\sigma$ to denote the CDF of $\mathbb{P}_{h, n} * \mathcal{N}(0, \sigma^2), \mathbb{P}_{h} * \mathcal{N}(0, \sigma^2), \mathcal{N}(0, \sigma^2)$. We will prove the results for $\epsilon$ sufficiently small, since cases of larger $\epsilon$ are direct corollary of cases of small $\epsilon$. 

	We fix $\sigma, K$, and let $\delta = \delta(\sigma, K, \epsilon)$ such that
	\begin{equation}\label{eq-definition-delta}\frac{(1 + \delta)(1 + \sigma^{2}/K^{2})^{2}}{2(1 - \delta)(1 + \sigma^{2}/K^{2}) - 4\delta\sigma^{2}/K^{2}} = \frac{(1 + \sigma^{2}/K^{2})^{2}}{2 + 2\sigma^{2}/K^{2}} + 2\epsilon,\end{equation}
	and we let $\zeta = \frac{\left(\frac{1}{2} + \frac{\sigma^2}{2K^2}\right)^2}{2\sigma^2}$. Then we know that $\lim_{\epsilon\to 0}\delta(\sigma, K, \epsilon) = 0$, and $\zeta - \frac{1}{2K^{2}} = \frac{\left(\frac{1}{2} - \frac{\sigma^2}{2K^2}\right)^2}{2\sigma^2} > 0$.
	With loss of generality we assume $\delta < \frac{1}{2}$. Therefore, for sufficiently small $\epsilon > 0$, we will have $\delta = \delta(\sigma, K, \epsilon) < \min\left\{\frac{1}{2}, 1 - \frac{1}{2K^{2}\zeta}\right\}$.
	\par We first show that for sufficiently large $h$ and some specific choice of $t\in (0, h-2)$, we will have
	\begin{align*}
		\mathbf{P}(X\in [t, t+2]) & \le\frac{4}{\sqrt{2\pi}\sigma}\exp\left(-(1-\delta)\zeta h^2\right)\quad \text{and}\quad \mathbf{P}(X\in [t+1, t+2]) \ge\frac{1}{2\sqrt{2\pi}\sigma}\exp\left(-(1+\delta)\zeta h^2\right).
	\end{align*}
	Actually, we have the following estimation of the probability of $\mathbb{P}_{h} * \mathcal{N}(0, \sigma^2)$ within the intervals $[t, t+2]$ and $[t+1, t+2]$: for $X\sim \mathbb{P}_h * \mathcal{N}(0, \sigma^2)$ and $t\in (0, h-2)$, we have
    \begin{align*}
        \mathbf{P}(X\in [t, t+2]) & \le 2\cdot \max_{t'\in [t, t+2]}\left[ \frac{1-p_h}{\sqrt{2\pi}\sigma}\exp\left(-\frac{t'^2}{2\sigma^2}\right) + \frac{p_h}{\sqrt{2\pi}\sigma}\exp\left(-\frac{(h - t')^2}{2\sigma^2}\right)\right]\\
        & \le \frac{2}{\sqrt{2\pi}\sigma}\cdot \left[\exp\left(-\frac{t^2}{2\sigma^2}\right) + \exp\left(-\frac{h^2}{2K^2} - \frac{(h - t - 2)^2}{2\sigma^2}\right)\right],\\
        \mathbf{P}(X\in [t + 1, t+2]) & \ge \min_{t'\in [t + 1, t+2]}\left[ \frac{1-p_h}{\sqrt{2\pi}\sigma}\exp\left(-\frac{t'^2}{2\sigma^2}\right) + \frac{p_h}{\sqrt{2\pi}\sigma}\exp\left(-\frac{(h - t')^2}{2\sigma^2}\right)\right] \ge \frac{1}{2\sqrt{2\pi}\sigma}\exp\left(-\frac{(t + 2)^2}{2\sigma^2}\right),
    \end{align*}
    where we have use the fact that $1 - p_{h}\ge \frac{1}{2}$ for all $h\ge 2K$. 
    Next, we would like to select the value of $t$ such that $-\frac{h^2}{2K^2} - \frac{(h - t)^2}{2\sigma^2} = -\frac{t^{2}}{2\sigma^{2}}$, e.g. $t = \frac{h}{2} + \frac{\sigma^2h}{2K^2}$. Since $\sigma < K$, we notice that $\exists \bar{h} > 0$ depending on $\sigma, K$ such that for $h > \bar{h}$, we have $t\in (0, h-2)$,
    and when $h$ goes to infinity, both $t$ and $h - t$ go to infinity as well. Hence for any $0 < \delta < 1$ there exists $C_1, C_h$ only depending on $K, \sigma$ and $\delta$ such that when $h > C_{h}$, we have $\frac{(h - t - 2)^2}{2\sigma^2}\le \frac{(1-\delta)(h - t)^{2}}{2\sigma^{2}}$ and $\frac{(t+2)^{2}}{2\sigma^{2}}\le \frac{(1 + \delta)t^{2}}{2\sigma^{2}}$, which implies that 
    \begin{align*}
        \mathbf{P}(X\in [t, t+2])\le \frac{4}{\sqrt{2\pi}\sigma}\cdot \exp\left(-\frac{(1-\delta)t^2}{2\sigma^2}\right) = \frac{4}{\sqrt{2\pi}\sigma}\exp\left(-\frac{(1-\delta)\left(\frac{1}{2} + \frac{\sigma^2}{2K^2}\right)^2h^2}{2\sigma^2}\right),\\
        \mathbf{P}(X\in [t+1, t+2])\ge \frac{1}{2\sqrt{2\pi}\sigma}\cdot \exp\left(-\frac{(1+\delta)t^2}{2\sigma^2}\right) = \frac{1}{2\sqrt{2\pi}\sigma}\exp\left(-\frac{(1+\delta)\left(\frac{1}{2} + \frac{\sigma^2}{2K^2}\right)^2h^2}{2\sigma^2}\right)
    \end{align*}
    holds for all $h > C_h$. We let $C_{1}\triangleq \frac{4}{\sqrt{2\pi}\sigma}$. 
    
    \par We next show that fix $h > 0$, for any $n$ such that $np_{h}\ge 128$ and $0 < t < h-2$, we have with probability at least $\frac{1}{16}$,
    $$\tilde{F}_{n, \sigma}(t) - F_\sigma(t)\ge \frac{1}{\sqrt{18n}}\exp\left(-\frac{h^2}{4K^2}\right).$$
    We first notice that for $0 < t < h$, $\tilde{F}_{n, \sigma}(t) - F_\sigma(t) = (\hat{p}_{h} - p_{h})(\Phi_\sigma(t - h) - \Phi_\sigma(t))$. Letting $U_i = \mathbf{1}_{X_k = h}$, according to Berry-Esseen Theorem~\citep{berry1941accuracy, esseen1956moment, durrett2019probability}, for $V\sim \mathcal{N}(0, 1)$, we have 
$$\sup_{x}\left|\mathbf{P}\left(\frac{1}{\sqrt{n\text{Var}[U_{1}]}}\sum_{l=1}^{n}[U_{l} - \mathbb{E}U_{1}]\le - x\right) - \mathbf{P}(V\le - x)\right|\le \frac{\mathbb{E}|U_{1} - \mathbb{E}[U_{1}]|^{3}}{2\sqrt{n}\sqrt{\text{Var}[U_{1}]}^{3}}.$$
When $p_{h} < 1/2$, we have $\mathbb{E}[U_1] = p_{h}$, $\text{Var}[U_{1}] = p_{h}(1 - p_{h})\ge p_{h}/2$ and $\mathbb{E}|U_{1} - \mathbb{E}[U_{1}]|^{3} \le \mathbb{E}|U_{1}|^{3} = \mathbb{E}[U_{1}] = p_{h}$. We choose $x = 1$, and noticing that $P(V > 1)\ge \frac{1}{8}$ we obtain
$$\mathbf{P}\left(\hat{p}_{h} - p_{h}\le -\sqrt{\frac{p_{h}}{2n}}\right) = \mathbf{P}\left(\frac{1}{n}\sum_{l=1}^{n}U_{l} - \mathbb{E}[U_{1}]\le -\sqrt{\frac{p}{2n}}\right)\ge \frac{1}{8} - \frac{\mathbb{E}|U_{1} - \mathbb{E}[U_{1}]|^{3}}{2\sqrt{n}\sqrt{\text{Var}[U_{1}]}^{3}}\ge \frac{1}{8} - \frac{1}{\sqrt{2np_{h}}}.$$
This indicates that $\hat{p}_{h} - p_{h}\le - \frac{1}{\sqrt{2n}}\exp\left(-h^2/(4K^2)\right)$ holds with probability at least $\frac{1}{8} - \frac{1}{\sqrt{2np_{h}}}.$ Then due to the fact that when $0 < t < h-2$ and $h > \sigma$, $\Phi_\sigma(t - h) - \Phi_\sigma(t)\le \Phi_\sigma(0) - \Phi_\sigma(h)\le -\frac{1}{3}$, we have with probability at least $\frac{1}{8} - \frac{1}{\sqrt{2np_{h}}}$, $\tilde{F}_{n, \sigma}(t) - F_\sigma(t)\ge 1/\sqrt{18n}\exp\left(-h^2/(4K^2)\right).$ Therefore when $np_{h}\ge 128$ we will have the above inequality holds with probability at least $\frac{1}{16}$. Combining these two results above, we know that whenever $h > \max\{C_h, \bar{h}, \sigma\}$ and $n$ satisfies that
    \begin{equation}\label{eq:pr2_a}
    	np_{h}\ge 128, \qquad \frac{1}{\sqrt{18n}}\exp\left(-\frac{h^2}{4K^2}\right) = \frac{4}{\sqrt{2\pi}\sigma}\exp\left(-(1-\delta)\zeta h^2\right),\end{equation}
    we will have for $t = \frac{h}{2} + \frac{\sigma^2h}{2K^2}$,
    $$\tilde{F}_{n, \sigma}(t) - F_\sigma(t)\ge \frac{1}{\sqrt{18n}}\exp\left(-\frac{h^2}{4K^2}\right)\ge \frac{4}{\sqrt{2\pi}\sigma}\exp\left(-(1-\delta)\zeta h^2\right)\ge \mathbf{P}(X\in [t, t+2])$$
    and hence $\tilde{F}_{n, \sigma}(t)\ge F_\sigma(t+2)$ holds with probability at least $1/16.$ Hence Lemma~\ref{lem2} indicates that with probability at least $1/16$ we have
    $$W_2(\mathbb{P}_{h}*\mathcal{N}(0, \sigma^{2}), \mathbb{P}_{h, n_{h}}*\mathcal{N}(0, \sigma^{2}))^2\ge \mathbf{P}(X\in [t+1, t+2])\ge \frac{1}{2\sqrt{2\pi}\sigma}\exp\left(-(1+\delta)\zeta h^2\right)$$
    and hence
    \begin{equation}\label{eq-bound-w}
    	\mathbb{E}\left[W_2(\mathbb{P}_{h}*\mathcal{N}(0, \sigma^{2}), \mathbb{P}_{h, n_{h}}*\mathcal{N}(0, \sigma^{2}))\right]\ge \frac{1}{16\sqrt{2\sqrt{2\pi}\sigma}}\exp\left(-\frac{(1+\delta)\zeta h^2}{2}\right) = C_{1}n^{ -\alpha - \epsilon},\end{equation}
    where $C_{1} = C_{1}(K, \sigma, \delta)$ is a positive constant. Here we use the fact that $\zeta = \left(\frac{1}{2} + \frac{\sigma^2}{2K^2}\right)^2/(2\sigma^2)\ge \frac{(\sigma/K)^{2}}{2\sigma^{2}} = \frac{1}{2K^{2}}$, then the second equation in \eqref{eq:pr2_a} indicates that
    \begin{equation}\label{eq-definition-h}
    	h = \sqrt{\left((1 - \delta)\zeta - \frac{1}{4K^{2}}\right)^{-1}\log\frac{12\sqrt{n}}{\sqrt{\pi}\sigma}}\end{equation}
    and when $\delta < \frac{1}{2}$ and $n\ge \pi\sigma^{2}/144$ this is well-defined. Bringing in this formula of $n$ into $\exp\left(-\frac{(1 + \delta)\zeta h^{2}}{2}\right)$ will result in the RHS in \eqref{eq-bound-w}, after noticing the definition of $\delta$ in \eqref{eq-definition-delta}.
    
    \par Finally we show that for sufficiently large $n$, there always exists an $h$ such that both~\eqref{eq:pr2_a} and also $h > \max\{C_h, \bar{h}, \sigma\}$ holds. For $n\ge \sqrt{\pi}\sigma/12$, we choose $h$ in \eqref{eq-definition-h} and the second equation in~\eqref{eq:pr2_a} holds, and also when $\delta < \frac{1}{2}$ there exists $n_{0}$ such that for any $n > n_{0}$ we have $h > \max\{C_{h}, \bar{h}, \sigma\}$. With this choice of $h$, we further have $np_{h} = n\exp\left(-h^{2}/(2K^{2})\right) = \Theta\left(n^{1 - (4K^{2}(1-\delta)\zeta - 1)^{-1}}\right),$ and since $\delta$ satisfies that $\delta < 1 - \frac{1}{2\zeta K^{2}}$ we know that $1 - (4K^{2}(1-\delta)\zeta - 1)^{-1} > 0$. Hence there exists a threshold $n_{0}$ such that for any $n > n_{1}$, we will have $np_{h}\ge 128$. Therefore, when $n > \max\{n_{0}, n_{1}\}$, with the choice of $h$ in \eqref{eq-definition-h}, we will have both~\eqref{eq:pr2_a} and also $h > \max\{C_h, \bar{h}, \sigma\}$ holds.
    
	\par Therefore, for any $\epsilon > 0$, we have $\mathbb{E}\left[W_2(\mathbb{P}_{h}*\mathcal{N}(0, \sigma^2), \mathbb{P}_{h, n}*\mathcal{N}(0, \sigma^2))\right] = \Omega\left(n^{-\alpha-\epsilon}\right)$.
\end{proof}

\subsection{LSI and $T_2$ constants for Bernoulli-Gaussian Mixtures}\label{apx:t2}

\subsubsection{Proof of the Non-Existence of Uniform Bound of LSI Constants for Bernoulli Distributions in \ref{sec-bernoulli}}
In this subsection, we will prove that for the Bernoulli distribution class in Section \ref{sec-bernoulli}, there constants in the corresponding log-Sobolev inequalities do not have a uniform bound. 
\begin{theorem}\label{thm-LSI}
Suppose $\sigma$ is a given constant which is smaller than $K$. Consider the following Bernoulli distributions:
$$\mathbb{P}_{h} = (1 - p_{h})\delta_{0} + p_{h}\delta_{h},\quad p_{h} = \exp\left(-h^{2}/(2K^{2})\right).$$
We use $C_{h}$ to denote the constant of LSI of distribution $\mu_{h} = \mathbb{P}_{h} * \mathcal{N}(0, \sigma^{2})$: $C_{h}$ is the smallest constant such that for any smoothed, compact supported function $f$ such that $\int_{\mathbb{R}}f^{2}d\mu_{h} = 1$, we have $\int_{\mathbb{R}}f^{2}\log f^{2}d\mu_{h}\le C_{h}\int_{\mathbb{R}}|f'|^{2}d\mu$. Then we have $\sup_{h\in\mathbb{R}_{+}}C_{h} = \infty.$
\end{theorem}
\begin{proof}[Proof of Theorem \ref{thm-LSI}]
We choose $x_{1} < -1 < 0 < x_{2} < h-1$, where $x_{1}$ and $x_{2}$ are determined later, and we let
$$f_{h}(x) = \begin{cases}
	0 & \quad x\le x_{1},\\
	t(x - x_{1}) & \quad x_{1}\le x\le x_{1} + 1,\\
	t &\quad x_{1} + 1\le x\le x_{2},\\
	-t(x - x_{2} - 1) &\quad x\ge x_{2},
\end{cases}$$
where $t$ is the constant chosen such that $\int_{\mathbb{R}}f_{h}^{2}d\mu_{h} = 1$. Then $f_{h}$ is a continuous function on $\mathbb{R}$, and $|f_{h}'(x)|\le t$ for any $x\in\mathbb{R}$. (Notice here $f_{h}$ is not a smooth function, but it has only finite points which are not smoothed. Hence after some simple smoothing procedure near these points, e.g. convolved with some mollifier, we can construct a sequence of functions converging to $f_{h}$ such that if the LSI works for functions in this sequence, the LSI also works for $f_{h}$.) Next, we will calculate the lower bound of $C_{h}$ such that the LSI works for function $f_{h}$. We denote
\begin{align*}
	& q_{h, 1} = \mu_{h}((-\infty, x_{1}]),\quad q_{h, 2} = \mu_{h}((x_{1}, x_{1} + 1]),\quad q_{h, 3} = \mu_{h}((x_{1}+1, x_{2}]),\\
	& q_{h, 4} = \mu_{h}((x_{2}, x_{2} + 1]),\quad q_{h, 5} = \mu_{h}((x_{2} + 1, \infty)).\end{align*}
Then we have $q_{h, 1} + q_{h, 2} + q_{h, 3} + q_{h, 4} + q_{h, 5} = 1$. According to the definition of $f$, we have $(q_{h, 2} + q_{h, 3} + q_{h, 4})t^{2}\ge \int_{\mathbb{R}}f_{h}^{2}d\mu_{h} = 1$, which indicates that $t^{2}\ge \frac{1}{q_{h, 2} + q_{h, 3} + q_{h, 4}}\ge 1$. Since for any $a\ge 0$, we have $a\log a\ge -1$, we also have
$$\int_{\mathbb{R}}f_{h}^{2}\log f_{h}^{2}d\mu_{h}\ge q_{h, 3}t^{2}\log t^{2} - (q_{h, 2} + q_{h, 4})\ge f_{h}^{2}d\mu_{h}\ge q_{h, 3}t^{2}\log t^{2} - (q_{h, 2} + q_{h, 4})t^{2}.$$
Moreover, we also notice that $|f_{h}'(x)|^{2} = t^{2}$ if $x\in (x_{1}, x_{1} + 1)\cup (x_{2}, x_{2} + 1)$, while $|f_{h}'(x)|^{2} = 0$ for other $x$. Therefore, we obtain that $\int_{\mathbb{R}}|f_{h}'|^{2}d\mu_{h} = (q_{h, 2} + q_{h, 4})t^{2}$. Hence if we require the LSI with constant $C_{h}$ holds for $f_{h}$, we will have
$$q_{h, 3}t^{2}\log t^{2} - (q_{h, 2} + q_{h, 4})t^{2}\le C_{h}(q_{h, 2} + q_{h, 4})t^{2},$$
which indicates that 
\begin{align*}
	C_{h} & \ge \frac{q_{h, 3}\log t^{2}}{q_{h, 2} + q_{h, 4}} - 1\ge \frac{- q_{h, 3}\log(q_{h, 2} + q_{h, 3} + q_{h, 4})}{q_{h, 2} + q_{h, 4}} - 1 = \frac{- q_{h, 3}\log(1 - q_{h, 1} - q_{h, 5})}{q_{h, 2} + q_{h, 4}} - 1\ge \frac{q_{h, 3}q_{h, 5}}{q_{h, 2} + q_{h, 4}} - 1.
\end{align*}
We use $\varphi_{\sigma^2}(x)$ to denote the PDF of $\mathcal{N}(0, \sigma^{2})$ at point $x$. According to the definition of $\mu_{h}$, and also noticing that $0 < x_{1} < h - 1$, we have $q_{h, 4} = \int_{x_{1}}^{x_{1} + 1}(1 - p_{h})\varphi_{\sigma^2}(x) + p_{h}\varphi_{\sigma^2}(x - h)dx\le \varphi_{\sigma^2}(x) + p_{h}\varphi_{\sigma^2}(h - x - 1)$, and also $q_{h, 5} = \int_{x_{1} + 1}^{\infty}(1 - p_{h})\varphi_{\sigma^2}(x) + p_{h}\varphi_{\sigma^2}(x - h)dx\ge \int_{x_{1} + 1}^{\infty}p_{h}\varphi_{\sigma^2}(x - h)dx\ge \int_{h}^{\infty}p_{h}\varphi_{\sigma^2}(x - h)dx = \frac{p_{h}}{2}$. We further notice that $\lim_{x_{1}\to-\infty}q_{h, 1} = \lim_{x_{1}\to -\infty}q_{h, 2} = 0.$ Hence letting $x_{1}\to -\infty$, we will obtain that $C_{h}$ satisfies
$$C_{h}\ge \lim_{x_{1}\to-\infty}\frac{q_{h, 3}q_{h, 5}}{q_{h, 2} + q_{h, 4}} - 1 = \lim_{x_{1}\to - \infty}\frac{q_{3}q_{5}}{q_{4}} - 1 = \frac{(1 - q_{4} - q_{5})q_{5}}{q_{4}} - 1\ge \frac{(1 - q_{5})q_{5}}{q_{4}} - 2.$$
When $\sigma < K$, we will choose $x = h\sqrt{\sigma/K}$, then we will have $\lim_{h\to\infty} x - h - 1 = \infty$, which indicates that
$$0\le \lim_{h\to\infty}\frac{q_{h, 4}}{p_{h}} = \lim_{h\to\infty}\frac{\varphi_{\sigma^2}(h\sqrt{\sigma/K}) + p_{h}\exp\varphi(h(1 - \sqrt{\sigma/K}))}{p_{h}} = 0,$$
and also $0\le \lim_{h\to\infty}q_{h, 5}\le \lim_{h\to\infty}\int_{h\sqrt{\sigma/K} + 1}^{\infty}\varphi_{\sigma^2}(x)dx + \lim_{h\to\infty}p_{h} = 0$, which indicates that $\lim_{h\to\infty}(1 - q_{h, 5}) = 1$. Above all, we obtain that $\lim_{h\to\infty}\frac{(1 - q_{5})q_{5}}{q_{4}} - 2 = \infty$, which indicates that $\lim_{h\to\infty}C_{h} = \infty$, and the uniform bound for $C_{h}$ does not exists.  
\end{proof}

\subsubsection{Proof of the Transportation-Entropy Inequality Constant}
\begin{theorem}\label{thm-transportation}
Suppose $\sigma$ is a given constant which is smaller than $K$. Consider the following Bernoulli distributions:
$$\mathbb{P}_{h} = (1 - p_{h})\delta_{0} + p_{h}\delta_{h},\quad p_{h} = \exp\left(-\frac{h^{2}}{2K^{2}}\right).$$
We use $C_{h}'$ to denote the constant of transportation-entropy inequality : $C_{h}$ is the smallest constant such that
\begin{equation}\label{eq-transportation}
	W_{2}(\mathbb{P}_{h} * \mathcal{N}(0, \sigma^{2}), \mathbb{Q})\le C_{h}'D_{KL}(\mathbb{P}_{h} * \mathcal{N}(0, \sigma^{2})\|\mathbb{Q})\qquad \forall\ \text{distribution}\ \mathbb{Q}.\end{equation}
Then we have $\sup_{h\in\mathbb{R}_{+}}C_{h}' = \infty$.
\end{theorem}
\begin{proof}
	We let $\mathbb{Q}_{h} = (1 - q_{h})\delta_{0} + q_{h}\delta_{h}$ with $q_{h} = p_{h} - \exp\left(-\frac{(1-\delta)(1 + \sigma^{2}/K^{2})^{2}h^{2}}{8\sigma^{2}}\right)$ for some $\delta$ smaller enough such that $(1 - \delta)(1 + \sigma^{2}/K^{2})^{2}h^{2} > 4\sigma^{2}/K^{2}$. According to data-processing inequality we have
	\begin{align*}
		& \quad D_{KL}(\mathbb{P}_{h} * \mathcal{N}(0, \sigma^{2})\|\mathbb{Q}_{h} * \mathcal{N}(0, \sigma^{2}))\le D_{KL}(\mathbb{P}_{h}\|\mathbb{Q}_{h}) = p_{h}\log \frac{p_{h}}{q_{h}} + (1-p_{h})\log \frac{1-p_{h}}{1-q_{h}}\\
		& = -p_{h}\log\left(1 + \frac{q_{h} - p_{h}}{p_{h}}\right) - (1 - p_{h})\log\left(1 + \frac{p_{h} - q_{h}}{1 - p_{h}}\right)\\
		& \le - p_{h}\cdot \frac{q_{h} - p_{h}}{p_{h}} + p_{h}\cdot \frac{(q_{h} - p_{h})^{2}}{p_{h}^{2}} - (1 - p_{h})\cdot \frac{p_{h} - q_{h}}{1-p_{h}} + (1-p_{h})\cdot \frac{(q_{h} - p_{h})^{2}}{(1-p_{h})^{2}}\le 2\exp\left(\frac{h^{2}}{2K^{2}}\right)(p_{h} - q_{h})^{2},
	\end{align*}
	where in the second inequality we use the fact that $-\log(1+x)\le -x + x^{2}$ for $x\ge -1/2$ and $\frac{q_{h} - p_{h}}{p_{h}}\ge -1/2$. Similar to the proof of Proposition \ref{prop-bernoulli}, and noticing that $F_{q, h}(t) - F_{p, h}(t) = (q_h - p_h)(\Phi_{\sigma}(t) - \Phi_{\sigma}(t - h))$ where $F_{q, h}, F_{p, h}, \Phi_{\sigma}$ are CDFs of distribution $\mathbb{Q}_h * \mathcal{N}(0, \sigma^{2}), \mathbb{P}_h * \mathcal{N}(0, \sigma^{2}), \mathcal{N}(0, \sigma^{2})$. We can prove that 
	$$W_{2}(\mathbb{P}_h * \mathcal{N}(0, \sigma^{2}), \mathbb{Q}_h * \mathcal{N}(0, \sigma^{2}))^{2} = \Omega\left(\exp\left(-\frac{(1-\delta)(1 + \sigma^{2}/K^{2})^{2}h^{2}}{8\sigma^{2}}\right)\right)$$
	while
	$$D_{KL}(\mathbb{P}_{h} * \mathcal{N}(0, \sigma^{2})\| \mathbb{Q}_h * \mathcal{N}(0, \sigma^{2})) = \mathcal{O}\left(\frac{h^{2}}{2K^{2}} - \frac{(1-\delta)(1 + \sigma^{2}/K^{2})^{2}h^{2}}{4\sigma^{2}}\right).$$
	Since $(1 - \delta)(1 + \sigma^{2}/K^{2})^{2}h^{2} > 4\sigma^{2}/K^{2}$, letting $h\to\infty$ we obtain that $\sup_{h\in\mathbb{R}_{+}}C_{h}' = \infty.$
\end{proof}

\section{Missing Proofs}\label{sec: appendix-proofs}
\subsection{Proof of Proposition \ref{prop: lb-lb-ub}}
We let $X = Y + Z$, where $Y\in \mathbb{P}, Z\sim \mathcal{N}(0, \sigma^2)$ are independent. We first prove \eqref{eq: lb-lb} by noticing that
\begin{align*}
	&\hspace{-0.5cm} 2\sigma^2\log\mathbf{P}(X\in [t_kr_{k} + 1, t_kr_{k} + 2])\ge 2\sigma^2\log\mathbf{P}(Y = r_k, Z\in [(t_k - 1)r_k + 1, (t_k - 1)r_k + 2])\\
	& \ge 2\sigma^2\log p_{k} + 2\sigma^2\log \mathbf{P}(Z\in [(t_k-1)r_{k} + 1, (t_k-1)r_{k} + 2]) \stackrel{(i)}{=} -\frac{r_{k}^{2}}{2K^{2}} - \frac{(t_k-1)^{2}(r_{k} + 2)^{2}}{2\sigma^{2}} - C_l\\
	& \ge -\left(\kappa + (t_k-1)^{2}\right)(r_{k} + 2)^{2} - C_l = -\left(t_k^{2} - \kappa
	c_k - c_k\right) (r_{k} + 2)^{2} - C_l,
\end{align*}
for some constant $C_l$, where in $(i)$ we use the definition of $p_k$ in \eqref{eq: def-p-k}. 

Next we prove \eqref{eq: lb-ub}. We first notice that
\begin{align*}
	& \hspace{-0.5cm} \mathbf{P}(X\in [t_kr_{k}, t_kr_{k} + 2]) = \sum_{j=0}^{\infty} \mathbf{P}(Y = r_j, Z\in [t_kr_k - r_j, t_kr_k - r_j + 2]) = \sum_{j=0}^{\infty}p_{j}\cdot \mathbf{P}(Z\in [t_kr_{k} - r_{j}, t_kr_{k} - r_{j} + 2])\\
	& \lesssim \sum_{j=0}^{k}p_{j}\exp\left(-\frac{(t_kr_{k} - r_{j})^{2}}{2\sigma^{2}}\right) + \sum_{j=k+1}^{\infty}p_{j}\exp\left(-\frac{(t_kr_{k} - r_{j} + 2)^{2}}{2\sigma^{2}}\right)
\end{align*} 
Then we upper bound the above inequality for $j = 0$, $0 < j< k$, $j = k$ and $j > k$ separately:
\begin{enumerate}
	\item For $j = 0$, noticing that $t_k^2\ge t_k^2 - c_k\kappa - c_k$, we have
	$$p_{j}\exp\left(-\frac{(t_kr_{k} - r_{j})^{2}}{2\sigma^{2}}\right)\le \exp\left(-\frac{t_kr_{k}^{2}}{2\sigma^{2}}\right)\lesssim\exp\left(- \frac{(t_k^{2} - \kappa c_k - c_k)(r_{k} - 2)^{2}}{2\sigma^{2}}\right).$$
	\item For $j > k$, we have $r_{j}^{2}\ge r_{k+1}^{2} + j - (k+1)$. When $c_k\ge \frac{\kappa + 3}{1 - \kappa}$, we have $c_k - t_k\ge 1$, which implies
	\begin{align*}
		(t_kr_{k} - r_{j} + 2)^{2} \ge ((c_k - t_k)r_{k} - 2)^{2}\ge (c_k - t_k)^{2}(r_{k} - 2)^{2} = (t_k^{2} - c_k\kappa - c_k - c_k^{2}\kappa)(r_{k} - 2)^{2}.
	\end{align*}
	For $j > k$, we further have 
	$$r_j^2\ge j - (k+1) + r_{k+1}^2 = j - (k+1) + \frac{K^2 \kappa c_k^2r_k^2}{\sigma^2}\ge j - (k+1) + \frac{K^2 \kappa c_k^2(r_k - 2)^2}{\sigma^2}.$$ 
	Hence we obtain:
	\begin{align*}
		& \hspace{-0.5cm} \sum_{j=k+1}^{\infty}\exp\left(-\frac{r_{j}^{2}}{2K^{2}} - \frac{(t_jr_{k} - r_{j} + 2)^{2}}{2\sigma^{2}}\right) \le \left(\sum_{j=k+1}^{\infty}\exp\left(-\frac{j-(k+1)}{2K^{2}}\right)\right)\cdot \exp\left(-(t_k^{2} - c_k\kappa - c_k)\cdot\frac{(r_{k} - 2)^{2}}{2\sigma^{2}}\right)\\
		& \lesssim \exp\left(-(t_k^{2} - c_k\kappa - c_k)\cdot \frac{(r_{k} - 2)^{2}}{2\sigma^{2}}\right).
	\end{align*}
	\item For $j < k$, we first notice $(t_kr_{k} - r_{j})^{2}\ge \left(t_k - \frac{1}{c_{k-1}}\right)^2r_k^2\ge \left(t_k^2 - \frac{2t_k}{c_{k-1}}\right)r_k^2\ge \left(t_k^2 - \kappa c_k - c_k\right)r_k^2$, where the first inequality uses $r_{j}\le \frac{r_{k}}{c_{k-1}}$ and the last inequality uses the definition of $t_k$ in \eqref{eq: def-t-k}. Therefore we have
	\begin{align*}
		& \quad \sum_{j=1}^{k-1}\exp\left(-\frac{r_{j}^{2}}{2K^{2}} - \frac{(t_kr_{k} - r_{j})^{2}}{2\sigma^{2}}\right)\le \sum_{j=1}^{k-1}\exp\left(-\frac{j}{2K^{2}} - \frac{(t_k^{2} - \kappa c_k - c_k)r_{k}^{2}}{2\sigma^{2}}\right)\\
		& \lesssim\exp\left(- \frac{(t_k^{2} - \kappa c_k - c_k)r_{k}^{2}}{2\sigma^{2}}\right) \le \exp\left(- \frac{(t_k^{2} - \kappa c_k - c_k)(r_{k} - 2)^{2}}{2\sigma^{2}}\right).
	\end{align*}
	\item For $j = k$, with the choice of $t_k = \frac{1}{2}(1 + c_k)(1 + \kappa)$ in \eqref{eq: def-t-k}, we obtain:
	\begin{align*}
		\exp\left(-\frac{r_{k}^{2}}{2K^{2}} - \frac{(t_kr_{k} - r_{k})^{2}}{2\sigma^{2}}\right) \lesssim \exp\left(-\left(t_k^{2} - \kappa c_k - c_k\right)\cdot \frac{r_{k}^{2}}{2\sigma^{2}}\right)\le\exp\left(-\left(t_k^{2} - \kappa c_k - c_k\right)\cdot \frac{(r_{k} - 2)^{2}}{2\sigma^{2}}\right).	\end{align*}
\end{enumerate}
Therefore, we have
$$\mathbf{P}(X\in [t_kr_{k}, t_kr_{k} + 2])\lesssim \exp\left(-(t_k^{2} - c_k\kappa - c_k)\cdot \frac{(r_{k} - 2)^{2}}{2\sigma^{2}}\right),$$
which implies \eqref{eq: lb-ub}.

\end{document}